\documentstyle[12pt]{article}
\input amssym.def
\begin{document}
\def \Z{\Bbb Z}
\def \C{\Bbb C}
\def \R{\Bbb R}
\def \Q{\Bbb Q}
\def \N{\Bbb N}
\def \wt{{\rm wt}}
\def \tr{{\rm tr}}
\def \span{{\rm span}}
\def \Res{{\rm Res}}
\def \End{{\rm End}}
\def \Ind {{\rm Ind}}
\def \Irr {{\rm Irr}}
\def \Aut{{\rm Aut}}
\def \Hom{{\rm Hom}}
\def \mod{{\rm mod}}
\def \ann{{\rm Ann}}
\def \<{\langle} 
\def \>{\rangle} 
\def \a{\alpha }
\def \e{\epsilon }
\def \l{\lambda }
\def \L{\Lambda }
\def \g{\gamma}
\def \b{\beta }
\def \om{\omega }
\def \ch{\chi}
\def \s{\sigma}

\def \bconj{\begin{conj}\label}
\def \econj{\end{conj}}
\def \be{\begin{equation}\label}
\def \ee{\end{equation}}
\def \bex{\begin{exa}\label}
\def \eex{\end{exa}}
\def \bl{\begin{lem}\label}
\def \el{\end{lem}}
\def \bt{\begin{thm}\label}
\def \et{\end{thm}}
\def \bp{\begin{prop}\label}
\def \ep{\end{prop}}
\def \br{\begin{rem}\label}
\def \er{\end{rem}}
\def \bc{\begin{coro}\label}
\def \ec{\end{coro}}
\def \bd{\begin{de}\label}
\def \ed{\end{de}}
\def \pf{{\bf Proof. }}
\def \voa{{vertex operator algebra}}

\newtheorem{thm}{Theorem}[section]
\newtheorem{prop}[thm]{Proposition}
\newtheorem{coro}[thm]{Corollary}
\newtheorem{conj}[thm]{Conjecture}
\newtheorem{exa}[thm]{Example}
\newtheorem{lem}[thm]{Lemma}
\newtheorem{rem}[thm]{Remark}
\newtheorem{de}[thm]{Definition}
\newtheorem{hy}[thm]{Hypothesis}
\makeatletter
\@addtoreset{equation}{section}
\def\theequation{\thesection.\arabic{equation}}
\makeatother
\makeatletter

\newcommand{\rw}{\rightarrow}
\newcommand{\n}{\:^{\circ}_{\circ}\:}

\begin{center}{\Large \bf 
The regular representation, Zhu's $A(V)$-theory and induced modules}
\end{center}

\begin{center}{Haisheng Li\footnote{Supported in part by NSF grants 
DMS-9616630 and DMS-9970496}\\
Department of Mathematical Sciences\\
Rutgers University-Camden\\
Camden, NJ 08102}
\end{center}

\begin{abstract}
The regular representation is related to Zhu's $A(V)$-theory and
an induced module from an $A(V)$-module to a $V$-module is defined
in terms of the regular representation. As an application, 
a new proof of Frenkel and Zhu's fusion rule theorem is obtained.
\end{abstract}

\section{Introduction}

In a remarkable paper [Z],
 Zhu constructed  among other things an associative 
algebra $A(V)$ for each vertex operator algebra $V$ and established
a one-to-one correspondence between the set of equivalence classes of
irreducible $A(V)$-modules and the set of equivalence classes of
lowest weight irreducible generalized $V$-modules. With this one-to-one 
correspondence, the classification of irreducible
$V$-modules is reduced to the classification 
of irreducible $A(V)$-modules. 
In [FZ], Zhu's $A(V)$-theory was extended further to determine
fusion rules by using $A(V)$-modules and bimodules associated to
$V$-modules. Since Zhu had developed his $A(V)$-theory, 
there have been many applications and generalizations 
(see for examples [A1-2], [DLM1-4], [DMZ], [DN1-3], [FZ], [KW], [W]). 
In Zhu's one-to-one correspondence, the functor from a weak $V$-module 
to an $A(V)$-module is a restriction with respect to both the space 
and the algebra,
and the functor from an $A(V)$-module to a (weak) $V$-module is,
to a certain extent, analogous to the induction functor in group theory.

In Lie group theory, for a Lie group $G$ and a subgroup $H$,
the induced $G$-module from an $H$-module $U$ is defined (cf. [Ki])
to be
$$\Ind_{H}^{G}U
=\{ f : G \rightarrow U\;|\; f(hg)=hf(g)\;\;\mbox{ for }h\in H,\; g\in G\},$$
where $(gf)(g')=f(g'g)$ for $g,g'\in G,\; f\in \Ind_{G}^{H}U$.
The construction of the induced module can be explained as follows:
First, $L^{2}(G)$ or $C^{0}(G)$ is (naturally) a $G\times G$-module.
(Certain $G\times G$-submodules are the modules affording the regular 
representation of $G$.)
More generally, for any (finite-dimensional) vector space $U$, the space
 $C^{o}(G,U)$ of continuous functions from $G$ to $U$ 
is a $G\times G$-module.
Second, the subspace $\Ind_{H}^{G}(U)$ of (left) $H$-invariant
functions from $G$ to $U$ is a $G$-submodule of 
$C^{0}(G,U)$ viewed as
a $G$-module through the identification $G=G\times 1$.

In [Li3] we defined regular representations 
of vertex operator algebras and established certain results.
More specifically, for a vertex operator algebra $V$ and a nonzero 
complex number $z$, we constructed a (weak) $V\otimes V$-module
${\cal{D}}_{P(z)}(V)$ out of the full dual space $V^{*}$ of $V$, and
we obtained certain results of Peter-Weyl type.
Note that unlike in group theory, there is no natural $V\otimes V$-module 
structure on $V^{*}$.
In view of this, ${\cal{D}}_{P(z)}(V)$ 
in a sense plays the role of $C^{0}(G)$. 

The main purpose of this paper is to relate Zhu's $A(V)$-theory 
to the regular representation in the spirit
of the induced module theory for a Lie group.
First, for a vector space $U$,
we construct a (weak) $V\otimes V$-module ${\cal{D}}_{P(z)}(V,U)$,
a subspace of $\Hom(V,U)$, which plays the role of $C^{0}(G,U)$.
Note that in Zhu's $A(V)$-theory, $A(V)$
is not a subalgebra of $V$ in the usual sense and
$A(V)$ does not naturally act on the whole space of a (weak) 
$V$-module.
In view of this, for an $A(V)$-module $U$, it does not make sense
to consider $A(V)$-invariant functions from $V$ to $U$.
On the other hand, given a (weak) $V$-module $W$, 
there is a canonical $A(V)$-bimodule
$A(W)$ [FZ], constructed as a quotient space of $W$ just as $A(V)$ 
is a quotient space of $V$ (see Section 3 for the definition);
and there is an $A(V)$-module $\Omega(W)$, a subspace of $W$. 
By definition, $\Omega(W)$ consists of those $w$ 
such that $v_{n}w=0$ for homogeneous $v\in V$ and for $n\ge \wt v$.
(Of course, $\Omega(W)$ can also be considered as the invariant space
with respect to a certain Lie algebra.)
In the case that $W$ is a lowest weight irreducible generalized $V$-module,
$\Omega(W)$ is the lowest weight subspace. 

We here define an induced module using the following
restriction-expansion strategy.
Since $A(V)$ is a quotient space of $V$, 
any linear function from $A(V)$ to $U$ lifts to
a linear function from $V$ to $U$. Then
we first restrict ourselves to linear functions from $V$ to $U$, 
which are lifted from linear functions from $A(V)$ to $U$, 
or simply just linear functions from $A(V)$ to $U$.
Now, it makes perfect sense to consider
(left) $A(V)$-invariant functions from $A(V)$ to $U$.
It is a classical fact that the space $\Hom(A(V),U)$ of linear 
functions from $A(V)$ to $U$ is a natural $A(V)$-module
containing the space $\Hom_{A(V)}(A(V),U)$ of $A(V)$-invariant
linear functions from $A(V)$ to $U$ as a submodule.
Of course,  $\Hom_{A(V)}(A(V),U)$ is 
canonically isomorphic to $U$.
On the other hand, it is shown 
(Proposition \ref{ptopWU}, Theorem \ref{tembeding}) that 
$\Hom(A(V),U)$ is a subspace of
${\cal{D}}_{P(-1)}(V,U)$, moreover
$\Hom(A(V),U)$ and 
$\Omega({\cal{D}}_{P(-1)}(V,U))$ $(\subset \Hom(V,U))$
coincide as natural $A(V)\otimes A(V)$-modules. 
To summarize, we have the following information:
\begin{eqnarray}
U=\Hom_{A(V)}(A(V),U)\subset \Hom(A(V),U)
=\Omega({\cal{D}}_{P(-1)}(V,U)).
\end{eqnarray}
Then we define the induced module $\Ind_{A(V)}^{V}U$ 
to be the submodule of ${\cal{D}}_{P(-1)}(V,U)$
generated by $\Hom_{A(V)}(A(V),U)$ $(=U)$ under the action of 
$V\otimes {\C}$.

Note that the results of [Li3] were more general than what we needed
for regular representations. For any weak $V$-module $W$,
a weak $V\otimes V$-module  ${\cal{D}}_{P(z)}(W)$
was constructed and it was proved that
the fusion rule of
type ${W'\choose W_{1}W_{2}}$ is equal to
$$\dim \Hom_{V\otimes V}(W_{1}\otimes W_{2},
{\cal{D}}_{P(-1)}(W))$$
for generalized $V$-modules $W,W_{1},W_{2}$.
Furthermore, if $W_{1}$ and $W_{2}$ are
lowest weight generalized $V$-modules, 
it was shown (Corollary 4.6, [Li3]) that 
the fusion rule of
type ${W'\choose W_{1}W_{2}}$ is equal to
$$\dim \Hom_{A(V)\otimes A(V)}(W_{1}(0)\otimes W_{2}(0), 
\Omega({\cal{D}}_{P(-1)}(W))),$$
where $W_{1}(0)$ and $W_{2}(0)$ are the corresponding 
lowest weight subspaces.
It is proved (Proposition \ref{ptopWU} and Theorem \ref{tembeding}) that 
$\Omega({\cal{D}}_{P(-1)}(W))$ and $A(W)^{*}$
coincide as natural $A(V)\otimes A(V)$-modules
for any weak $V$-module $W$. 
Using these results, we obtain a new proof 
of Frenkel and Zhu's fusion rule theorem\footnote{The original
theorem [FZ] was corrected in [Li1-2] (see
Corollary \ref{cintertwining} below).} 
which 
asserts that the fusion rule of type ${W_{2}\choose WW_{1}}$
for irreducible $V$-modules $W, W_{1}, W_{2}$ is equal to
$$\dim \Hom_{A(V)}(A(W)\otimes_{A(V)}W_{1}(0),W_{2}(0))$$
under a certain condition (Corollary \ref{cintertwining}).

In [DLin], an induced module theory for a vertex operator algebra 
with respect to a vertex operator subalgebra 
was established.
Let $V_{1}$ be a vertex operator subalgebra of $V$
and let $U$ be an irreducible  $V_{1}$-module. 
In general, $U$ could lift to either a $V$-module or a so-called
twisted $V$-module by an automorphism of $V$, but not both.
(A $V$-module is a twisted module corresponding to the identity
automorphism.) In this regard,
this theory is quite different from and more complicated than
the classical theory. We hope to study Dong-Lin's
induced module theory in terms of regular representations later.

We would like to thank Professor James Lepowsky for 
providing me many valuable suggestions,
as he has generously done for many of my papers. We also would like
to thank Professor Geoffrey Mason (the editor) 
for providing a list of corrections on the exposition.

The paper is organized as follows: In Section 2, we review
the construction of the weak $V\otimes V$-module ${\cal{D}}_{P(z)}(W)$
and the main results, and then construct a weak $V\otimes V$-module
${\cal{D}}_{P(z)}(W,U)$. In Section 3, we 
identify $\Hom(A(W),U)$ with $\Omega({\cal{D}}_{P(z)}(W,U))$
as natural $A(V)\otimes A(V)$-modules, and we define the
induced $V$-module $\Ind_{A(V)}^{V}U$ for a given $A(V)$-module $U$.
In Section 4, we give a new proof of the Frenkel and Zhu's fusion
rule theorem. 

\section{Weak $V\otimes V$-modules ${\cal{D}}_{P(z)}(W)$ 
and ${\cal{D}}_{P(z)}(W,U)$}

In this section we shall first review the construction of the weak 
$V\otimes V$-module ${\cal{D}}_{P(z)}(W)$ and the main results from [Li3],
and then construct a (weak) $V\otimes V$-module
${\cal{D}}_{P(z)}(W,U)$ as a generalization.

We use standard definitions and notations as given in [FLM] and [FHL].
A vertex operator algebra is denoted by $V$, or by $(V,Y,{\bf 1},\omega)$
with more information,
where ${\bf 1}$ is the vacuum vector and $\omega$ is the Virasoro element.
We also use the notion of weak module as defined in [DLM2]---A weak module
satisfies all the axioms given in [FLM] and [FHL] for the notion 
of a module except that no grading is required.

We typically use letters $x,y, x_{1},x_{2},\dots$ for mutually commuting
formal variables and $z,z_{0},\dots$ for complex numbers.
For a vector space $U$, $U[[x,x^{-1}]]$ is the vector space of all
(doubly infinite) formal series with coefficients in $U$ and
$U((x))$ is the space of formal Laurent series. Sometimes we also
use $U[x,x^{-1}]]$ for $U((x^{-1}))$.
We emphasize the following standard formal variable convention:
\begin{eqnarray}
& &(x_{1}-x_{2})^{n}=\sum_{i\ge 0}(-1)^{i}{n\choose i}x_{1}^{n-i}x_{2}^{i},\\
& &(x-z)^{n}=\sum_{i\ge 0}(-z)^{i}{n\choose i}x^{n-i},\\
& &(z-x)^{n}=\sum_{i\ge 0}(-1)^{i}z^{n-i}{n\choose i}x^{i}
\end{eqnarray}
for $n\in {\Z},\; z\in {\C}^{\times}$.

Recall the following simple result from [Li3]:

\bl{lbasic5}
Let $U$ be a vector space, $U_{1}$ a subspace and let
\begin{eqnarray}
f(x)=\sum_{n\in {\Z}}f_{n}x^{-n-1}\in U[[x,x^{-1}]],
\;\;g(x)=\sum_{n\in {\Z}}g_{n}x^{-n-1}
\in U_{1}[[x,x^{-1}]].
\end{eqnarray}
Suppose that either $f(x)\in U((x))$ or $f(x)\in U((x^{-1}))$ and that
there exist $k\in {\N}$ and $z\in {\C}^{\times}$ such that
\begin{eqnarray}\label{ebasic52}
(x-z)^{k}f(x)=(x-z)^{k}g(x).
\end{eqnarray}
Then for $n\in {\Z}$,
\begin{eqnarray}
f_{n}\in \mbox{\rm linear span }\{g_{m}\;|\;m\ge n\}
\end{eqnarray}
if $f(x)\in U((x))$ and
\begin{eqnarray}
f_{n}\in \mbox{\rm linear span }\{g_{m}\;|\;m\le n\}
\end{eqnarray}
if $f(x)\in U((x^{-1}))$. In particular, $f(x)\in U_{1}[[x,x^{-1}]]$.
\el

For vector spaces $U_{1},U_{2}$, a linear map $f\in \Hom(U_{1},U_{2})$ 
extends canonically to a linear map from $U_{1}[[x,x^{-1}]]$ to 
$U_{2}[[x,x^{-1}]]$. We shall use this canonical extension
without any comments.

Let $V$ be a vertex operator algebra. 
For $v\in V$, we set (cf. [FHL], [HL1])
\begin{eqnarray}
Y^{o}(v,x)=Y(e^{xL(1)}(-x^{-2})^{L(0)}v,x^{-1}).
\end{eqnarray}
For a weak $V$-module $W$, $Y^{o}(v,x)$ lies in $\Hom (W,W[x,x^{-1}]])$ 
because $e^{xL(1)}(-x^{-2})^{L(0)}v\in V[x,x^{-1}]$ and
$Y(u,x^{-1})w\in W[x,x^{-1}]]$ for $u\in V,\; w\in W$. 
More generally,  for any complex number $z_{0}$,
$Y^{o}(v,x+z_{0})$  lies in $\Hom (W,W[x,x^{-1}]])$, where by definition
\begin{eqnarray}
Y^{o}(v,x+z_{0})w=(Y^{o}(v,y)w)|_{y=x+z_{0}}
\end{eqnarray}
for $w\in W$.
Let $W$ be a weak $V$-module and let $U$ be a vector space, e.g., $U={\C}$.
For $v\in V,\;f\in \Hom(W,U)$, the compositions
$fY^{o}(v,x)$ and $fY^{o}(v,x+z_{0})$ for any complex number $z_{0}$
are elements of $(\Hom (W,U))[[x,x^{-1}]]$. 

Now let us review the main definitions and results about 
${\cal{D}}_{P(z)}(W)$ from [Li3]. 

\bd{drecall1}
{\em  [Li3] Let $V$ be a vertex operator algebra, $W$ a weak $V$-module
and $z$ a nonzero complex number. Define ${\cal{D}}_{P(z)}(W)$
to be the subspace of $W^{*}$, consisting of those $\alpha$ such that
for each $v\in V$, there exist $k,l\in {\N}$ such that for $w\in W$,
\begin{eqnarray}
x^{l}(x-z)^{k}\<\alpha, Y^{o}(v,x)w\>\in {\C}[x],
\end{eqnarray}
 or what is equivalent, the series
$\<\alpha, Y^{o}(v,x)w\>$, an element of ${\C}[x,x^{-1}]]$,
 absolutely converges
in the domain $|x|>|z|$ to a rational function of the form 
$x^{-l}(x-z)^{-k}g(x)$, where $g(x)\in {\C}[x]$.}
\ed

The following is an obvious characterization 
for $\alpha$ lying in ${\cal{D}}_{P(z)}(W)$ 
without involving matrix-coefficients.

\bl{lrecall1} {\rm [Li3]}
Let $W, z$ be given as before and let $\alpha\in W^{*}$. 
Then $\alpha\in {\cal{D}}_{P(z)}(W)$ if and only if
for $v\in V$, there exist $k,l\in {\N}$ such that
\begin{eqnarray}\label{ers1}
x^{l}(x-z)^{k}\alpha Y^{o}(v,x)\in W^{*}[[x]],
\end{eqnarray}
or equivalently, if and only if for $v\in V$, 
there exists $k\in {\N}$ such that
\begin{eqnarray}
(x-z)^{k}\alpha Y^{o}(v,x)\in W^{*}((x)).
\end{eqnarray}
\el

Let ${\C}(x)$ be the algebra of rational functions of $x$.
The $\iota$-maps $\iota_{x;0}$ and $\iota_{x;\infty}$ from 
${\C}(x)$ to ${\C}[[x,x^{-1}]]$ are defined as follows:
for any rational function $f(x)$, 
$\iota_{x;0}f(x)$ is the Laurent series expansion of $f(x)$ at $x=0$
and $\iota_{x;\infty}f(x)$ is the Laurent series expansion 
of $f(x)$ at $x=\infty$. These are injective 
${\C}[x,x^{-1}]$-linear maps. In terms of
the formal variable convention, we have
\begin{eqnarray}
& &\iota_{x;0}\left((x-z)^{n}f(x)\right)=(-z+x)^{n}\iota_{x;0}f(x),\\
& &\iota_{x;\infty}\left((x-z)^{n}f(x)\right)=(x-z)^{n}\iota_{x;\infty}f(x)
\end{eqnarray}
for $n\in {\Z},\; z\in {\C}^{\times},\; f(x)\in {\C}(x)$.

{}From the definition, 
for $\alpha\in {\cal{D}}_{P(z)}(W),\; v\in V,\; w\in W$,
$\<\alpha, Y^{o}(v,x)w\>$ lies in the range of $\iota_{x;\infty}$.
Then $\iota_{x;\infty}^{-1}\<\alpha, Y^{o}(v,x)w\>$ is a well 
defined element of ${\C}(x)$.

\bd{drecall2} {\rm [Li3]}
{\em For $v\in V$, $\alpha\in {\cal{D}}_{P(z)}(W)$, we define
$$Y_{P(z)}^{L}(v,x)\alpha,\;\;\;\; 
Y_{P(z)}^{R}(v,x)\alpha\in W^{*}[[x,x^{-1}]]$$
by
\begin{eqnarray}
\<Y_{P(z)}^{L}(v,x)\alpha,w\>&=&\iota_{x;0}\left(\iota_{x;\infty}^{-1}
\<\alpha,Y^{o}(v,x+z)w\>\right)\\
\<Y_{P(z)}^{R}(v,x)\alpha,w\>&=&\iota_{x;0}\iota_{x;\infty}^{-1}
\<\alpha,Y^{o}(v,x)w\>
\end{eqnarray}
for $w\in W$.}
\ed

\bl{lrecall1'} {\rm [Li3]}
Let $v\in V,\; \alpha\in {\cal{D}}_{P(z)}(W)$. Then
\begin{eqnarray}
& &(-z+x)^{k}Y_{P(z)}^{R}(v,x)\alpha=(x-z)^{k}\alpha Y^{o}(v,x),\\
& &(z+x)^{l}Y_{P(z)}^{L}(v,x)\alpha=(x+z)^{l}\alpha Y^{o}(v,x+z),
\end{eqnarray}
where $k$ and $l$ are any pair of nonnegative integers such that (\ref{ers1})
 holds.
\el

We have ([Li3], Proposition 3.24):

\bp{precall2'} {\rm [Li3]} 
Let $W$ be a weak $V$-module and let $z$ be a nonzero complex number.
Then 
\begin{eqnarray}
Y_{P(z)}^{L}(v,x)\alpha,\;\;\;\; 
Y_{P(z)}^{R}(v,x)\alpha\in ({\cal{D}}_{P(z)}(W))((x))
\end{eqnarray}
for $v\in V,\; \alpha \in {\cal{D}}_{P(z)}(W)$. Furthermore,
\begin{eqnarray}
Y_{P(z)}^{L}(u,x_{1})Y_{P(z)}^{R}(v,x_{2})
=Y_{P(z)}^{R}(v,x_{2})Y_{P(z)}^{L}(u,x_{1})
\end{eqnarray}
on ${\cal{D}}_{P(z)}(W)$ for $u,v\in V$.
\ep

In view of Proposition \ref{precall2'}, $Y^{L}_{P(z)}$ and $Y^{R}_{P(z)}$
give rise to a well defined linear map 
\begin{eqnarray}
Y_{P(z)}=Y^{L}_{P(z)}\otimes Y^{R}_{P(z)}: V\otimes V\rightarrow 
\left(\End\;{\cal{D}}_{P(z)}(W)\right)[[x,x^{-1}]].
\end{eqnarray}
Then we have ([Li3], Theorem 3.17, 
Propositions 3.21 and 3.24 and Theorem 3.25):

\bt{trecall2} {\rm [Li3]} 
Let $W$ be a weak $V$-module and let $z$ be a nonzero complex number.
Then the pairs $({\cal{D}}_{P(z)}(W), Y_{P(z)}^{L})$ and
$({\cal{D}}_{P(z)}(W), Y_{P(z)}^{R})$ carry 
the structure of a weak $V$-module and
the pair $({\cal{D}}_{P(z)}(W), Y_{P(z)})$ 
carries the structure of a weak $V\otimes V$-module.
\et

For a ${\C}$-graded vector space $M=\coprod_{h\in {\C}}M_{(h)}$,
following [HL1] we define the formal completion
\begin{eqnarray}
\overline{M}=\prod_{h\in {\C}}M_{(h)}.
\end{eqnarray}
Recall from [FHL] that $M'=\coprod_{h\in {\C}}M_{(h)}^{*}$. Then 
\begin{eqnarray}
\overline{M'}=M^{*}.
\end{eqnarray}

We shall need the following notions.
A {\em generalized} $V$-module [HL1] is a weak $V$-module on which
$L(0)$ semisimply acts. 
 Then for a generalized $V$-module $W$ we have 
the $L(0)$-eigenspace decomposition: $W=\coprod_{h\in {\C}}W_{(h)}$.
Thus, a generalized $V$-module 
satisfies all the axioms defining the notion of a $V$-module
([FLM], [FHL]) except the two grading restrictions on the 
homogeneous subspaces.
If a generalized $V$-module furthermore
satisfies the lower truncation condition (one of the two grading 
restrictions), we call it a {\em lower truncated} generalized module [H1].

Following [HL1], we choose a
branch $\log z$ of the log function so that
\begin{eqnarray}
\log z =\log |z|+i \arg z\;\;\;\mbox{ with }\;\;0\le \arg z< 2\pi,
\end{eqnarray}
and arbitrary values of the log function will be denoted by
\begin{eqnarray}
l_{p}(z)=\log z+2p\pi i
\end{eqnarray}
for $p\in {\Z}$.

Let  $W, W_{1}$ and $W_{2}$ be generalized $V$-modules and
let ${\cal{Y}}$ be an intertwining operator of type 
${W'\choose W_{1}W_{2}}$.
For $w_{(1)}\in W_{1},\; w_{(2)}\in W_{2}$, we set [HL1]
\begin{eqnarray}
{\cal{Y}}(w_{(1)},e^{l_{p}(z)})w_{(2)}
=\left({\cal{Y}}(w_{(1)},x)w_{(2)}\right)|_{x^{h}=e^{h l_{p}(z)}, h\in {\C}}
\in \overline{W'}\;(=W^{*}).
\end{eqnarray}
Note that  ${\cal{Y}}(w_{(1)},x)w_{(2)}$ 
in general involves non-integral, even complex powers of $x$.
We have ([Li3], Theorem 4.5):

\bp{precall4'} {\rm [Li3]}
Let  $W, W_{1}$ and $W_{2}$ be generalized $V$-modules,
${\cal{Y}}$ an intertwining operator of type ${W'\choose W_{1}W_{2}}$
and let $p\in {\Z}$. Then
\begin{eqnarray}
{\cal{Y}}(w_{(1)},e^{l_{p}(z)})w_{(2)}\in {\cal{D}}_{P(z)}(W)
\end{eqnarray}
for $w_{(1)}\in W_{1},\; w_{(2)}\in W_{2}$.
\ep

In view of Proposition \ref{precall4'}, for an intertwining operator 
${\cal{Y}}$ of type ${W'\choose W_{1}W_{2}}$ 
we have a linear map
\begin{eqnarray}
F_{{\cal{Y}},p}^{P(z)}: 
& &W_{1}\otimes W_{2}\rightarrow {\cal{D}}_{P(z)}(W)\nonumber\\
& &(w_{(1)},w_{(2)})\mapsto F_{{\cal{Y}},p}^{P(z)}(w_{(1)}\otimes w_{(2)})
={\cal{Y}}(w_{(1)},e^{l_{p}(z)})w_{(2)}
\end{eqnarray}
for $w_{(1)}\in W_{1},\; w_{(2)}\in W_{2}$.

For generalized $V$-modules $W, W_{1}$ and $W_{2}$, 
following [HL1] we denote by
${\cal{V}}^{W'}_{W_{1}W_{2}}$
the space of intertwining operators of type ${W'\choose W_{1}W_{2}}$.
Then we have ([Li3], Corollary 4.6):

\bt{trecall4} {\rm [Li3]}
Let $W, W_{1}$ and $W_{2}$ be lower truncated generalized $V$-modules, 
let $z$ be 
a nonzero complex number and let $p\in {\Z}$. Then the linear map
\begin{eqnarray}
F_{p}[P(z)]_{W_{1}W_{2}}^{W'}:& & {\cal{V}}^{W'}_{W_{1}W_{2}}\rightarrow 
\Hom_{V\otimes V} (W_{1}\otimes W_{2},{\cal{D}}_{P(z)}(W))\nonumber\\
& &{\cal{Y}}\mapsto F_{{\cal{Y}},p}^{P(z)}
\end{eqnarray}
is a linear isomorphism.
\et

Next, we shall generalize the notion of ${\cal{D}}_{P(z)}(W)$ 
by incorporating a vector space $U$.

\bd{dDWU}
{\em Let $W$ be a weak $V$-module, $U$ a vector space and 
$z$ a nonzero complex number.
Define ${\cal{D}}_{P(z)}(W,U)$ to be the subset of $\Hom(W,U)$,
consisting of each $f$ such that
for $v\in V$, there exist $k,l\in {\N}$ 
such that
\begin{eqnarray}\label{eDWUcharc}
(x-z)^{k}x^{l}\<u^{*}, fY^{o}(v,x)w\>\in {\C}[x]
\end{eqnarray}
for all $u^{*}\in U^{*},\; w\in W$, or what is equivalent,
for all $u^{*}\in U^{*},\; w\in W$, 
the formal series
$$\<u^{*}, fY^{o}(v,x)w\>,$$
an element of ${\C}[x,x^{-1}]]$,
absolutely converges in the domain $|x|>|z|$ to a rational function
of the form $x^{-l}(x-z)^{-k}g(x)$ for $g(x)\in {\C}[x]$.}
\ed

Clearly, ${\cal{D}}_{P(z)}(W,U)$ is a subspace of $\Hom(W,U)$.
When $U={\C}$, ${\cal{D}}_{P(z)}(W,{\C})$ gives us ${\cal{D}}_{P(z)}(W)$.

\bl{ldefdwuequiv}
Let $f\in \Hom(W,U)$. Then the following statements are equivalent:

(a) $f\in {\cal{D}}_{P(z)}(W,U)$.

(b) For $v\in V$, there exist $k,l\in {\N}$ such that
\begin{eqnarray}\label{edefDWU}
(x-z)^{k}x^{l}fY^{o}(v,x)\in (\Hom(W,U))[[x]].
\end{eqnarray}

(c)  For $v\in V$, there exist $k,l\in {\N}$ such that for each $w\in W$,
\begin{eqnarray}\label{edefDWU2}
(x-z)^{k}x^{l}fY^{o}(v,x)w\in U[x].
\end{eqnarray}
\el

\pf Clearly, (a) implies (b), and (c) implies (a). 
Since $Y^{o}(v,x)w\in W[x,x^{-1}]]$ for $v\in V,\;w\in W$,
we see that (b) implies (c). $\;\;\;\;\Box$

Let $v\in V,\; f\in {\cal{D}}_{P(z)}(W,U)$ and let $k,l\in {\N}$
be such that (\ref{edefDWU2}) holds.
Then by changing variable we get
\begin{eqnarray}\label{edefDWU3}
x^{k}(x+z)^{l}fY^{o}(v,x+z)w\in U[x]
\end{eqnarray}
for $w\in W$.

\bd{dLRactions}
{\em Let $W, U$ and $z$ be given as before. 
For $v\in V,\; f\in {\cal{D}}_{P(z)}(W,U)$, we define
two elements $Y^{L}_{P(z)}(v,x)f$ and $Y_{P(z)}^{R}(v,x)$
of $(\Hom (W,U))[[x,x^{-1}]]$ by
\begin{eqnarray}
(Y_{P(z)}^{L}(v,x)f)(w)
&=&(z+x)^{-l}\left((x+z)^{l}f(Y^{o}(v,x+z)w)\right)\\
(Y_{P(z)}^{R}(v,x)f)(w)&=&(-z+x)^{-k}\left((x-z)^{k}f(Y^{o}(v,x)w)\right)
\end{eqnarray}
for $w\in W$, where $k,l$ are any pair of (possibly negative) integers 
such that (\ref{edefDWU2}) holds.}
\ed

First, in view of (\ref{edefDWU2}) and (\ref{edefDWU3}), 
both $(z+x)^{-l}\left((x+z)^{l}f(Y^{o}(v,x+z)w)\right)$
and $(-z+x)^{-k}\left((x-z)^{k}f(Y^{o}(v,x)w)\right)$ lie
in $U((x))$, so that 
$Y^{L}_{P(z)}(v,x)f$ and $Y^{R}_{P(z)}(v,x)f$ make sense.
However, we are not allowed to
remove the left-right brackets to cancel $(x-z)^{k}$ or $(x+z)^{l}$ 
because of the nonexistence of terms $(z+x)^{-l}f(Y^{o}(v,x+z)w)$ and 
$(-z+x)^{-k}f(Y^{o}(v,x)w)$.
Second, they are also well defined, i.e., they are independent
of the choice of the pair of integers $k,l$. Indeed, 
if $k',l'$ are another pair of integers such that (\ref{edefDWU2}) holds,
say for example, $k\ge k'$, then 
\begin{eqnarray}
& &(-z+x)^{-k}\left((x-z)^{k}fY^{o}(v,x)w\right)\nonumber\\
&=&(-z+x)^{-k}\left((x-z)^{k-k'}(x-z)^{k'}fY^{o}(v,x)w\right)\nonumber\\
&=&(-z+x)^{-k}(x-z)^{k-k'}\left((x-z)^{k'}fY^{o}(v,x)w\right)\nonumber\\
&=&(-z+x)^{-k'}\left((x-z)^{k'}fY^{o}(v,x)w\right).
\end{eqnarray}

{}From definition we immediately have:

\bl{lDWUconn}
For $v\in V,\; f\in {\cal{D}}_{P(z)}(W,U)$,
\begin{eqnarray}
& &(z+x)^{l}Y_{P(z)}^{L}(v,x)f=(x+z)^{l}fY^{o}(v,x+z),\\
& &(-z+x)^{k}Y_{P(z)}^{R}(v,x)f=(x-z)^{k}fY^{o}(v,x),
\end{eqnarray}
where $k,l$ are any pair of integers 
such that (\ref{edefDWU2}) holds. $\;\;\;\;\Box$
\el

In terms of rational functions and the $\iota$-maps we immediately have
(cf. [DL], [FHL]):

\bl{lDWUconn2}
For $v\in V,\; f\in {\cal{D}}_{P(z)}(W,U),\; u^{*}\in U^{*},\; w\in W$,
\begin{eqnarray}
& &\<u^{*},(Y_{P(z)}^{L}(v,x)f)(w)\>
=\iota_{x;0}\iota_{x;\infty}^{-1}\<u^{*},fY^{o}(v,x+z)w\>,\\
& &\<u^{*},(Y_{P(z)}^{R}(v,x)f)(w)\>
=\iota_{x;0}\iota_{x;\infty}^{-1}\<u^{*},fY^{o}(v,x)w\>.\;\;\;\;\Box
\end{eqnarray}
\el
Let $W, U$ and $z$ be given as before.
Consider $U^{*}\otimes W$ as a weak $V$-module 
with the action of $V$ on $W$. 
Then in view of Theorem  \ref{trecall2} we have a
weak $V\otimes V$-module ${\cal{D}}_{P(z)}(U^{*}\otimes W)$.
Let $\alpha \in (U^{*}\otimes W)^{*}$. Then 
$\alpha\in {\cal{D}}_{P(z)}(U^{*}\otimes W)$ if and only if
for $v\in V$, there exist $k,l\in {\N}$ such that
\begin{eqnarray}
(x-z)^{k}x^{l}\<\alpha, u^{*}\otimes Y^{o}(v,x)w\>\in {\C}[x]
\end{eqnarray}
for all $u^{*}\in U^{*},\; w\in W$.

Let $\eta$ be the canonical embedding of
$\Hom(W,U)$ into $(U^{*}\otimes W)^{*}$, i.e.,
for $f\in \Hom(W,U)$, $u^{*}\in U^{*},\; w\in W$,
\begin{eqnarray}
\<\eta(f),u^{*}\otimes w\>=\<u^{*},f(w)\>.
\end{eqnarray}
Let $f\in {\cal{D}}_{P(z)}(W,U)\;(\subset \Hom(W,U))$.
For $v\in V$, let $l,k\in {\N}$ such that
\begin{eqnarray}
(x-z)^{k}x^{l}\<u^{*}, fY^{o}(v,x)w\>\in {\C}[x]
\end{eqnarray}
for all $u^{*}\in U^{*},\; w\in W$, that is,
\begin{eqnarray}
(x-z)^{k}x^{l}\<\eta(f), u^{*}\otimes Y^{o}(v,x)w\>\in {\C}[x]
\end{eqnarray}
for all $u^{*}\in U^{*},\; w\in W$. Then
$\eta(f)\in {\cal{D}}_{P(z)}(U^{*}\otimes W)$.
This proves
\begin{eqnarray}
\eta({\cal{D}}_{P(z)}(W,U))\subset {\cal{D}}_{P(z)}(U^{*}\otimes W).
\end{eqnarray}

On the other hand, let $f\in \Hom (W,U)$.
If $\eta(f)\in {\cal{D}}_{P(z)}(U^{*}\otimes W)$,
for $v\in V$, there exist $k,l\in {\N}$ such that
\begin{eqnarray}
(x-z)^{k}x^{l}\<\eta(f), u^{*}\otimes Y^{o}(v,x)w\>\in {\C}[x]
\end{eqnarray}
for all $u^{*}\in U^{*},\; w\in W$. That is,
\begin{eqnarray}
(x-z)^{k}x^{l}\<u^{*},f(Y^{o}(v,x)w)\>\in {\C}[x].
\end{eqnarray}
Then $f\in {\cal{D}}_{P(z)}(W,U)$. This shows
\begin{eqnarray}
\eta(\Hom(W,U))\cap {\cal{D}}_{P(z)}(U^{*}\otimes W)
\subset \eta({\cal{D}}_{P(z)}(W,U)).
\end{eqnarray}
Therefore, we have proved:

\bl{lDWU} Let $W$ be a weak $V$-module, $U$ a vector space
and $z$ a nonzero complex number. Then
\begin{eqnarray}
\eta({\cal{D}}_{P(z)}(W,U))
=\eta(\Hom(W,U))\cap {\cal{D}}_{P(z)}(U^{*}\otimes W).\;\;\;\;\Box
\end{eqnarray}
\el

Furthermore, we have:

\bp{pWU}
Let $W$ be a weak $V$-module, and let $U$ be a vector space. Then
$\eta({\cal{D}}_{P(z)}(W,U))$ is a weak $V\otimes V$-submodule of 
${\cal{D}}_{P(z)}(U^{*}\otimes W)$. Furthermore,
\begin{eqnarray}
\eta(Y^{L}_{P(z)}(v,x)f)=Y^{L}_{P(z)}(v,x)\eta(f),\;\;
\eta(Y_{P(z)}^{R}(v,x)f)=Y_{P(z)}^{R}(v,x)\eta(f)
\end{eqnarray}
for $v\in V,\;f\in {\cal{D}}_{P(z)}(W,U)$.
\ep

\pf Let $v\in V,\; f\in {\cal{D}}_{P(z)}(W,U)$. 
Since $\eta(f)\in {\cal{D}}_{P(z)}(U^{*}\otimes W)$ (Lemma \ref{lDWU}),
by Lemma \ref{lrecall1'} there exist $k,l\in {\N}$ such that
 \begin{eqnarray}
(x-z)^{k}Y^{R}_{P(z)}(v,x)\eta(f)
&=&(x-z)^{k}\eta(f)Y^{o}(v,x)\label{erelationyy*1}\\
(x+z)^{l}Y^{L}_{P(z)}(v,x)\eta(f)&=&(x+z)^{l}\eta(f)Y^{o}(v,x+z).
\label{erelationyy*2}
\end{eqnarray}
For $u^{*}\in U^{*},\; w\in W$,  we have
\begin{eqnarray}
\<\eta(f)Y^{o}(v,x), u^{*}\otimes w\>
&=&\<\eta(f),Y^{o}(v,x)(u^{*}\otimes w)\>\nonumber\\
&=&\<\eta(f),u^{*}\otimes Y^{o}(v,x)w\>\nonumber\\
&=&\<u^{*},fY^{o}(v,x)w\>\nonumber\\
&=&\<\eta(fY^{o}(v,x)),u^{*}\otimes w\>.
\end{eqnarray}
Then
\begin{eqnarray}\label{e2.55}
\eta(f)Y^{o}(v,x)=\eta(fY^{o}(v,x)).
\end{eqnarray}
Consequently,
\begin{eqnarray}
\eta(f)Y^{o}(v,x)\;\left(=\eta(fY^{o}(v,x))\right)
\in \eta(\Hom(W,U))[[x,x^{-1}]].
\end{eqnarray}
Then it follows from Lemma \ref{lbasic5} and 
(\ref{erelationyy*1})-(\ref{erelationyy*2}) that 
$$Y^{R}_{P(z)}(v,x)\eta(f),\;\;Y^{L}_{P(z)}(v,x)\eta(f)
\in \eta(\Hom(W,U))[[x,x^{-1}]],$$
so that from Lemma \ref{lDWU},
\begin{eqnarray}
Y^{R}_{P(z)}(v,x)\eta(f),\;\;Y^{L}_{P(z)}(v,x)\eta(f)
\in \eta({\cal{D}}_{P(z)}(W,U))[[x,x^{-1}]].
\end{eqnarray}
This proves that $\eta({\cal{D}}_{P(z)}(W,U))$ is a weak 
$V\otimes V$-submodule of ${\cal{D}}_{P(z)}(U^{*}\otimes W)$.

Let $v\in V,\;f\in {\cal{D}}_{P(z)}(W,U),\; u^{*}\in U^{*},\; w\in W$.
Then using Lemma \ref{lDWUconn2} we get
\begin{eqnarray}
\<Y^{L}_{P(z)}(v,x)\eta(f), u^{*}\otimes w\>
&=&\iota_{x;0}\iota_{x;\infty}^{-1}\<\eta(f),u^{*}\otimes Y^{o}(v,x+z)w)\>
\nonumber\\
&=&\iota_{x;0}\iota_{x;\infty}^{-1}\<u^{*}, f(Y^{o}(v,x+z)w)\>
\nonumber\\
&=&\<u^{*}, (Y^{L}_{P(z)}(v,x)f)w)\>\nonumber\\
&=&\<\eta(Y^{L}_{P(z)}(v,x)f), u^{*}\otimes w\>.
\end{eqnarray}
Thus
$$\eta(Y^{L}_{P(z)}(v,x)f)=Y^{L}_{P(z)}(v,x)\eta(f).$$
Similarly we can prove
$$\eta(Y^{R}_{P(z)}(v,x)f)=Y^{R}_{P(z)}(v,x)\eta(f).$$
This completes the proof.
$\;\;\;\;\Box$

In view of Theorem \ref{trecall2} and
Proposition \ref{pWU} we immediately have:

\bt{tDWU}
Let $W$ be a weak $V$-module, $U$ a vector space and $z$ a nonzero 
complex number. Then the pairs $({\cal{D}}_{P(z)}(W,U), Y_{P(z)}^{L})$ and
$({\cal{D}}_{P(z)}(W,U), Y_{P(z)}^{R})$ carry 
the structure of a weak $V\otimes V$-module and the actions
$Y_{P(z)}^{L}$ and $Y_{P(z)}^{R}$ of $V$ on ${\cal{D}}_{P(z)}(W,U)$
commute. Furthermore, set 
\begin{eqnarray}
Y_{P(z)}=Y_{P(z)}^{L}\otimes Y_{P(z)}^{R}.
\end{eqnarray}
Then the pair $({\cal{D}}_{P(z)}(W,U), Y_{P(z)})$ 
carries the structure of a weak $V\otimes V$-module. $\;\;\;\;\Box$
\et

In view of Proposition \ref{pWU} and (\ref{e2.55}), from 
([Li3], Proposition 3.22) 
we immediately have the following 
relations among $fY^{o}(v,x), Y^{L}(v,x)$ and $Y^{R}(v,x)f$:

\bc{cthreeconn}
Let $v\in V,\; f\in {\cal{D}}_{P(z)}(W,U)$. Then
\begin{eqnarray}\label{ethreeconn}
& &x_{0}^{-1}\delta\left(\frac{x-z}{x_{0}}\right)fY^{o}(v,x)
-x_{0}^{-1}\delta\left(\frac{z-x}{-x_{0}}\right)Y_{P(z)}^{R}(v,x)f\nonumber\\
&=&z^{-1}\delta\left(\frac{x-x_{0}}{z}\right)Y_{P(z)}^{L}(v,x_{0})f.
\end{eqnarray}
\ec

For convenience, from now on we shall drop
the ``$P(z)$'' from the notations $Y^{L}_{P(z)}$ and $Y^{R}_{P(z)}$
when there is no confusion.

\section{Zhu's $A(V)$-theory and induced module $\Ind_{A(V)}^{V}U$}
In this section, given a weak $V$-module $W$ and 
a nonzero complex number $z$, we construct an $A(V)\otimes A(V)$-module 
$A(W,z)$, generalizing Frenkel and Zhu's notion of $A(W)$,
and then we relate $\Hom(A(W,z),U)$ to a canonical subspace of
${\cal{D}}_{P(z)}(W,U)$. Using this connection, we define
the induced module $\Ind_{A(V)}^{V}U$ from an $A(V)$-module $U$.

First we define or review certain notions.
A {\em lowest weight} generalized $V$-module is a generalized $V$-module 
such that $W=\coprod_{n\in \N}W_{(h+n)}$ for some $h\in {\C}$ and
$W_{(h)}$ generates $W$.
Furthermore, if $W\ne 0$, we call the unique $h$ 
{\em the lowest weight} of $W$. 
An {\em ${\N}$-graded} weak $V$-module [Z] is a weak $V$-module $W$
together with an ${\N}$-grading  $W=\coprod_{n\in {\N}}W(n)$ 
such that
\begin{eqnarray}
v_{m}W(n)\subset W(n+\wt v-m-1)
\end{eqnarray}
for homogeneous $v\in V$ and for $m\in {\Z},\; n\in {\N}$, where
by definition $W(n)=0$ for $n<0$.
An {\em ${\N}$-gradable} weak $V$-module is a weak $V$-module $W$
on which there exists an ${\N}$-grading
such that $W$ together the grading becomes an ${\N}$-graded module.
A vertex operator algebra $V$ is said to be {\em rational} 
[Z] (cf. [DLM2])
if every ${\N}$-gradable weak $V$-module is a direct sum of
irreducible ${\N}$-gradable weak $V$-modules. There are also 
different definitions of rationality (see for example [HL1]).

Now we recall Zhu's construction of $A(V)$ and 
the main results from [Z]. 
Let $V$ be a vertex operator algebra. Set 
\begin{eqnarray}
O(V)=\mbox{linear span} \{ \Res_{x}x^{-2}(1+x)^{\wt u}Y(u,x)v
\;\mbox{ for homogeneous }u,v\in V\}.
\end{eqnarray}
Note that we {\em do not} assume that $V$ has the special property
that $V=\oplus_{n\ge 0}V_{(n)}$, so that $\wt u$ could be {\em negative}, 
hence the formal series $(1+x)^{\wt u}$ and 
$(x+1)^{\wt u}$ may be different.

For homogeneous $u,v\in V$, we define [Z]
\begin{eqnarray}
u*v=\Res_{x}x^{-1}(1+x)^{\wt u}Y(u,x)v
\left(=\sum_{i\ge 0}{\wt u\choose i}u_{i-1}v\right).
\end{eqnarray}
Then extend the definition of $*$ on $V$ by linearity. Set
\begin{eqnarray}
A(V)=V/O(V).
\end{eqnarray}
The following is 
the first of Zhu's theorems in his $A(V)$-theory.

\bp{pzhu1} {\rm [Z] }
Let $(V,Y,{\bf 1},\omega)$ be a vertex operator algebra. Then
the space $O(V)$ is a two-sided ideal of the nonassociative algebra
$(V,*)$ and the quotient algebra $A(V)$ ($=V/O(V)$) is an 
associative algebra
with ${\bf 1}+O(V)$ being the identity element and with $\omega+O(V)$ 
being a central element. Furthermore,
$A(V)$ has an involution (anti-automorphism) $\theta$ given by
\begin{eqnarray}
\theta (v)=e^{L(1)}(-1)^{L(0)}v\;\;\;\mbox{ for }v\in V.
\end{eqnarray}
\ep

Let $W$ be a weak $V$-module. Following [DLM2] we define
\begin{eqnarray}
\Omega(W)=\{ w\in W\;| \;v_{n}w=0\;\mbox{ for homogeneous }v\in V
\mbox{ and for }n\ge \wt v\}.
\end{eqnarray}
Equivalently, $w\in \Omega(W)$ if and only if $x^{\wt v}Y(v,x)w\in W[[x]]$
for each homogeneous $v\in V$.
Then we have ([Z],  [DLM2]):

\bp{pzhu2} For any weak $V$-module $W$, $\Omega(W)$ is a natural 
$A(V)$-module with $v+O(V)$ acting on $\Omega(W)$ as $v_{\wt v-1}$ 
for homogeneous $v\in V$. Furthermore, if $W=\coprod_{n\ge 0}W_{(n+h)}$
is a lowest weight irreducible generalized $V$-module with $W_{(h)}\ne 0$, 
then $\Omega(W)=W_{(h)}$ and it is an irreducible $A(V)$-module.
\ep

Let $W_{1}, W_{2}$ be weak $V$-modules and let $\psi$ be a 
$V$-homomorphism from $W_{1}$ to $W_{2}$. Clearly,
$\psi(\Omega(W_{1}))\subset \Omega(W_{2})$ and
the restriction $\Omega(\psi):=\psi|_{\Omega(W_{1})}$ is an 
$A(V)$-homomorphism. It is routine to check that
$\Omega$ is a functor from 
the category of weak $V$-modules to the category of $A(V)$-modules.
On the other hand, for any $A(V)$-module $U$ Zhu in [Z] 
constructed a ${\N}$-graded weak $V$-module $L(U)$ 
with $U=L(U)(0)\subset \Omega(L(U))$ (cf. [DLM2]).
Now we shall use the generalized regular representation of $V$ on 
${\cal{D}}_{P(z)}(V,U)$ to construct such an ${\N}$-graded weak
$V$-module.

Let $W$ be a weak $V$-module and let $z$ be a nonzero complex number.
Generalizing the definition of  $O(W)$ in [FZ], we define $O(W,z)$
to be the subspace of $W$, linearly spanned by elements
\begin{eqnarray}
\Res_{x}x^{-2}(1-zx)^{\wt v}Y(v,x)w
\end{eqnarray}
for homogeneous $v\in V$ and for $w\in W$. With this notion, $O(W)=O(W,-1)$.
Generalizing Frenkel and Zhu's left and right actions of $V$ on $W$ [FZ]
we define 
\begin{eqnarray}
v*_{P(z)}w&=&\Res_{x}(-z)^{-\wt v}x^{-1}(1-zx)^{\wt v}Y(v,x)w,
\label{epzavleft}\\
w*_{P(z)}v&=&\Res_{x}(-z)^{-\wt v}x^{-1}(1-zx)^{\wt v-1}Y(v,x)w
\label{epzavright}
\end{eqnarray}
for homogeneous $v\in V$ and for $w\in W$.
Then extend the definitions by linearity. (We recover
Frenkel and Zhu's actions when $z=-1$.) 
In the following we shall show that these generalized actions
actually are Frenkel and Zhu's actions of $V$ on $W$ with respect to
a {\em new module} structure.

\bl{lgeneralization}
Let $W$ be a weak $V$-module and let $z$ be a nonzero complex number.
For $v\in V$, set
\begin{eqnarray}
Y^{(z)}(v,x)=Y(z^{L(0)}v,zx).
\end{eqnarray}
Then $(W,Y^{(z)})$ carries the structure of a weak $V$-module.
Furthermore, for homogeneous $v\in V$ and for $m,n\in {\Z}$, we have
\begin{eqnarray}\label{eresidue}
\Res_{x}(-z)^{-\wt v}x^{m}(1-zx)^{n}Y(v,x)
=(-z)^{-m-1}\Res_{x}x^{m}(1+x)^{n}Y^{(-z^{-1})}(v,x).
\end{eqnarray}
\el

\pf From [FHL], we have
\begin{eqnarray}\label{eL(0)con}
z^{L(0)}Y(v,x)z^{-L(0)}=Y(z^{L(0)},zx)
\end{eqnarray}
on any generalized $V$-module (on which $L(0)$ semisimply acts).
In particular, this is true on the adjoint module $V$.
If $W$ is a generalized $V$-module, it follows 
immediately from (\ref{eL(0)con})
that $(W,Y^{(z)})$ is a weak $V$-module and it is isomorphic to $(W,Y)$
through the map $z^{L(0)}$.
For a general weak $V$-module $W$, replacing $(u,v)$ and 
$(x_{0},x_{1},x_{2})$
 by $(z^{L(0)}u,z^{L(0)}v)$ and $(zx_{0},zx_{1},zx_{2})$
in the Jacobi identity for $Y$, respectively, then
using (\ref{eL(0)con}) on $V$ we obtain
\begin{eqnarray}
& &z^{-1}x_{0}^{-1}\delta\left(\frac{x_{1}-x_{2}}{x_{0}}\right)
Y^{(z)}(u,x_{1})Y^{(z)}(v,x_{2})\nonumber\\
&-&z^{-1}x_{0}^{-1}\delta\left(\frac{x_{2}-x_{1}}{-x_{0}}\right)
Y^{(z)}(v,x_{2})Y^{(z)}(u,x_{1})\nonumber\\
&=&z^{-1}x_{2}^{-1}\delta\left(\frac{x_{1}-x_{0}}{x_{2}}\right)
Y(Y(z^{L(0)}u,zx_{0})z^{L(0)}v,zx_{2})\nonumber\\
&=&z^{-1}x_{2}^{-1}\delta\left(\frac{x_{1}-x_{0}}{x_{2}}\right)
Y(z^{L(0)}Y(u,x_{0})v,zx_{2})\nonumber\\
&=&z^{-1}x_{2}^{-1}\delta\left(\frac{x_{1}-x_{0}}{x_{2}}\right)
Y^{(z)}(Y(u,x_{0})v,x_{2}).
\end{eqnarray}
This proves the Jacobi identity for $Y^{(z)}$ while
the vacuum property and lower truncation condition clearly hold.
The identity (\ref{eresidue}) directly follows from changing variable
$y=-z^{-1}x$.$\;\;\;\;\Box$

With Lemma \ref{lgeneralization}, generalizations of
certain Zhu's theorems [Z], or Frenkel and Zhu's theorems [FZ]
will follow immediately. First, we have (cf. [Z]):

\bl{lzhuO(W)}
Let $W$ be a weak $V$-module and let $z$ be a nonzero complex number. Then
\begin{eqnarray}
\Res_{x}x^{-n-2}(1-zx)^{\wt v+m}Y(v,x)w\in O(W,z)
\end{eqnarray}
for homogeneous $v\in V$ and for $n\ge m\ge 0$, $w\in W$.$\;\;\;\;\Box$
\el

We also have:

\bl{lO(W)} Let $W$ and $z$ be given as before. Then
\begin{eqnarray}
\Res_{x}x^{-n-2}(1-zx)^{\wt v+m}Y(e^{x^{-1}L(1)}v,x)w\in O(W,z)
\end{eqnarray}
for any homogeneous $v\in V$ and for $n\ge m\ge 0,\; w\in W$.
\el

\pf Notice that  $\wt (L(1)^{i}v)=\wt v-i$ for $i\ge 0$.
Then using Lemma \ref{lzhuO(W)}, we get
\begin{eqnarray}
& &\Res_{x}x^{-n-2}(1-zx)^{\wt v+m}Y(e^{x^{-1}L(1)}v,x)w\nonumber\\
&=&\sum_{i\ge 0}\Res_{x}{1\over i!}x^{-n-i-2}(1-zx)^{\wt v+m}
Y(L(1)^{i}v,x)w
\nonumber\\
&=&\sum_{i\ge 0}\Res_{x}{1\over i!}x^{-n-i-2}(1-zx)^{(\wt v-i)+i+m}
Y(L(1)^{i}v,x)w\in O(W,z). \;\;\;\;\Box
\end{eqnarray}

Set $A(W,z)=W/O(W,z)$. Then $A(W)=A(W,-1)$. 
Noticing that from Lemma \ref{lgeneralization},
the left and right actions $*_{P(z)}$ are exactly the Frenkel and Zhu's
left and right actions on the module $(W,Y^{(-z^{-1})})$, we immediately
have:

\bp{p} {\rm [FZ]} Let $W$ be a weak $V$-module and let $z$ be 
a nonzero complex number. 
Then the left and right actions $*_{P(z)}$ of $V$ on $W$
defined in (\ref{epzavleft}) and (\ref{epzavright}) 
give rise to an $A(V)$-bimodule structure on $A(W,z)$.
\ep

We shall need the following result:

\bl{ltensor}
Let $V_{1}$ and $V_{2}$ be vertex operator algebras and let $E$
be a weak $V_{1}\otimes V_{2}$-module. Then
\begin{eqnarray}
\Omega_{V_{1}\otimes V_{2}}(E)=\Omega_{V_{1}}(E)\cap \Omega_{V_{2}}(E),
\end{eqnarray}
where $E$ is considered as a weak $V_{1}$-module and 
a weak $V_{2}$-module in the obvious way.
\el

\pf Clearly,
$$\Omega_{V_{1}\otimes V_{2}}(E)
\subset \Omega_{V_{1}}(E)\cap \Omega_{V_{2}}(E).$$
On the other hand, since the actions of $V_{1}$ and $V_{2}$ on
$E$ commute, 
\begin{eqnarray}\label{e3.18}
Y(v_{(2)},x)\Omega_{V_{1}}(E)\subset \Omega_{V_{1}}(E) [[x,x^{-1}]]
\end{eqnarray}
for $v_{(2)}\in V_{2}$.
Now, let $e\in \Omega_{V_{1}}(E)\cap \Omega_{V_{2}}(E)$
and let $v_{(1)}\in V_{1},\; v_{(2)}\in V_{2}$ be homogeneous.
Then
\begin{eqnarray}
x^{\wt v_{(2)}}Y(v_{(2)},x)e
\in E[[x]]\cap \Omega_{V_{1}}(E)[[x,x^{-1}]]=\Omega_{V_{1}}(E)[[x]],
\end{eqnarray}
so that using (\ref{e3.18}) we get
\begin{eqnarray}
x^{(\wt v_{(1)}\otimes v_{(2)})}Y(v_{(1)}\otimes v_{(2)},x)e
=x^{\wt v_{(1)}}Y(v_{(1)},x)\left(x^{\wt v_{(2)}}Y(v_{(2)},x)e\right)
\in E[[x]].
\end{eqnarray}
Thus $e\in \Omega_{V_{1}\otimes V_{2}}(E)$. This proves
$$\Omega_{V_{1}}(E)\cap \Omega_{V_{2}}(E)
\subset \Omega_{V_{1}\otimes V_{2}}(E)$$
and completes the proof. $\;\;\;\;\Box$

Now let $W$ be a weak $V$-module, $U$ a vector space and 
$z$ a nonzero complex number. 
Consider $\Hom(A(W,z),U)$ naturally as a subspace of $\Hom (W,U)$. 
Recall that $\eta$ is the canonical embedding of $\Hom (W,U)$
into $(U^{*}\otimes W)^{*}$. 
Then we have:

\bp{ptopWU}
Let $W$ be a weak $V$-module, $U$ a vector space, 
and $z$ a nonzero complex number. Then
\begin{eqnarray}
& &\Hom (A(W;z),U)\nonumber\\
&=&\Omega({\cal{D}}_{P(z)}(W,U))\\
&=&\{f\in \Hom(W,U)\;|\; 
x^{\wt v}(x-z)^{\wt v}f Y^{o}(v,x)\in (\Hom(W,U))[[x]]
\;\;\mbox{ for homogeneous }v\in V\}.\label{etopaw}
\end{eqnarray}
Furthermore, 
\begin{eqnarray}
& &(-z+x)^{\wt v}Y^{R}(v,x)f =(x-z)^{\wt v}f Y^{o}(v,x),\label{e3.23}\\
& &(z+x)^{\wt v}Y^{L}(v,x)f =(x+z)^{\wt v}f Y^{o}(v,x+z)\label{e3.24}
\end{eqnarray}
for $f\in \Hom (A(W,z),U)$ and for homogeneous $v\in V$.
\ep

\pf Let $T$ be the set defined in the right hand side of (\ref{etopaw}).
To prove the first assertion, in the following we shall prove
$$\Hom (A(W,z),U)\subset T\subset \Omega({\cal{D}}_{P(z)}(W,U))
\subset \Hom (A(W,z),U).$$
The second part follows immediately from Lemma \ref{lDWUconn}.

Let $f\in \Hom (A(W,z),U)\;(\subset \Hom(W,U))$ and 
let $v\in V$ be homogeneous. 
Then for $n\in {\N},\; w\in W$, 
by changing variable and using Lemma \ref{lO(W)} we get
\begin{eqnarray}
& &\Res_{x}x^{\wt v+n}(x-z)^{\wt v}fY^{o}(v,x)w\nonumber\\
&=&\Res_{x}x^{\wt v+n}(x-z)^{\wt v}
fY(e^{xL(1)}(-x^{-2})^{L(0)}v,x^{-1})w\nonumber\\
&=&\Res_{x}x^{-\wt v-n-2}(x^{-1}-z)^{\wt v}
fY(e^{x^{-1}L(1)}(-x^{2})^{L(0)}v,x)w\nonumber\\
&=&(-1)^{\wt v}f\left(\Res_{x}x^{-n-2}(1-zx)^{\wt v}
Y(e^{x^{-1}L(1)}v,x)w\right)\nonumber\\
&=&0.
\end{eqnarray}
This shows 
\begin{eqnarray}\label{echarc}
x^{\wt v}(x-z)^{\wt v}fY^{o}(v,x)\in (\Hom(W,U))[[x]].
\end{eqnarray}
That is, $f\in T$. Thus
$$\Hom (A(W,z),U)\subset T.$$
{}From the definition of ${\cal{D}}_{P(z)}(W,U)$, 
we immediately have
$$T\subset {\cal{D}}_{P(z)}(W,U).$$

Let $f\in T$ and let $v\in V$ be homogeneous.
By Lemma \ref{lDWUconn} we have
\begin{eqnarray}
& &(-z+x)^{\wt v}Y^{R}(v,x)f
=(x-z)^{\wt v}f Y^{o}(v,x),\\
& &(z+x)^{\wt v}Y^{L}(v,x)f
=(x+z)^{\wt v}f Y^{o}(v,x+z).
\end{eqnarray}
Then
\begin{eqnarray}
& &x^{\wt v}(-z+x)^{\wt v}Y^{R}(v,x)f
=x^{\wt v}(x-z)^{\wt v}f Y^{o}(v,x)
\in (\Hom(W,U))[[x]],\\
& &x^{\wt v}(z+x)^{\wt v}Y^{L}(v,x)f
=x^{\wt v}(x+z)^{\wt v}f Y^{o}(v,x+z)
\in (\Hom(W,U))[[x]].
\end{eqnarray}
We are also using (\ref{edefDWU3}). Hence
\begin{eqnarray}
& &x^{\wt v}Y^{R}(v,x)f
=(-z+x)^{-\wt v}\left[x^{\wt v}(-z+x)^{\wt v}Y^{R}(v,x)f\right]
\in (\Hom(W,U))[[x]],\\
& &x^{\wt v}Y^{L}(v,x)\alpha
=(z+x)^{-\wt v}\left[x^{\wt v}(z+x)^{\wt v}Y^{L}(v,x)f\right]
\in (\Hom(W,U))[[x]].
\end{eqnarray}
It follows from Lemma \ref{ltensor} that $f\in \Omega({\cal{D}}_{P(z)}(W,U))$. 
This proves
$$T \subset \Omega({\cal{D}}_{P(z)}(W,U)).$$

Let $f\in \Omega({\cal{D}}_{P(z)}(W,U))$ and let $v\in V$ be homogeneous.
Then
\begin{eqnarray}\label{ealphainomega}
x^{\wt v}Y^{L}(v,x)f, 
\;\;x^{\wt v}Y^{R}(v,x)f\in (\Hom(W,U))[[x]].
\end{eqnarray}
Multiplying (\ref{ethreeconn}) by $x^{\wt v}x_{0}^{\wt v}$, 
then taking $\Res_{x_{0}}$
(and using the fundamental properties of delta functions) we get
\begin{eqnarray}
& &x^{\wt v}(x-z)^{\wt v}fY^{o}(v,x)-x^{\wt v}(-z+x)^{\wt v}Y^{R}(v,x)f
\nonumber\\
&=&\Res_{x_{0}}z^{-1}\delta\left(\frac{x-x_{0}}{z}\right)
(z+x_{0})^{\wt v}x_{0}^{\wt v}Y^{L}(v,x_{0})f.
\end{eqnarray}
Then it follows from (\ref{ealphainomega}) that
\begin{eqnarray}
x^{\wt v}(x-z)^{\wt v}f Y^{o}(v,x)=x^{\wt v}(-z+x)^{\wt v}Y^{R}(v,x)f
\in (\Hom(W,U))[[x]].
\end{eqnarray}
That is, $f\in T$.
Furthermore,
for homogeneous $v\in V$ and for $w\in W$, since $(Y^{o})^{o}=Y$ [FHL],
we have
\begin{eqnarray}
& &\Res_{x}x^{-2}(1-zx)^{\wt v}fY(v,x)w\nonumber\\
&=&\Res_{x}x^{-2}(1-zx)^{\wt v}
fY^{o}(e^{xL(1)}(-x^{-2})^{L(0)}v,x^{-1})w\nonumber\\
&=&\Res_{x}(1-zx^{-1})^{\wt v}
fY^{o}(e^{x^{-1}L(1)}(-x^{2})^{L(0)}v,x)w\nonumber\\
&=&\Res_{x}(-1)^{\wt v}x^{\wt v}(x-z)^{\wt v}
fY^{o}(e^{x^{-1}L(1)}v,x)w\nonumber\\
&=&\sum_{i\ge 0}(-1)^{\wt v}{1\over i!}\Res_{x}x^{\wt v-i}(x-z)^{\wt v}
fY^{o}(L(1)^{i}v,x)w\nonumber\\
&=&\sum_{i\ge 0}(-1)^{\wt v}{1\over i!}
\Res_{x}x^{\wt (L(1)^{i}v)}(x-z)^{\wt (L(1)^{i}v)+i}
f Y^{o}(L(1)^{i}v,x)w\nonumber\\
&=&0
\end{eqnarray}
because
\begin{eqnarray}
& &\Res_{x}x^{\wt (L(1)^{i}v)}(x-z)^{\wt (L(1)^{i}v)+i}
f Y^{o}(L(1)^{i}v,x)w\nonumber\\
&=&\sum_{j=0}^{i}{i\choose j}
\Res_{x}x^{\wt (L(1)^{i}v)+j}(x-z)^{\wt (L(1)^{i}v)}
f Y^{o}(L(1)^{i}v,x)w\nonumber\\
&=&0.
\end{eqnarray}
This proves $f(O(W,z))=0$, hence $f\in \Hom(A(W,z),U)$.
Thus $\Omega({\cal{D}}_{P(z)}(W,U))\subset \Hom(A(W,z),U)$. This
completes the proof.$\;\;\;\;\Box$

It follows from Theorem \ref{tDWU}, Lemma \ref{ltensor}, and 
Proposition \ref{pzhu2}
that $\Omega({\cal{D}}_{P(z)}(W,U))$ is an $A(V)\otimes A(V)$-module.
On the other hand, because $A(W,z)$ is an $A(V)$-bimodule and $\theta$ is 
an involution of $A(V)$, from the classical fact
$\Hom (A(W,z),U)$ becomes an $A(V)\otimes A(V)$-module with
\begin{eqnarray}\label{etensormodule}
((a_{1},a_{2})f)(w)=f(\theta(a_{2})wa_{1})
\end{eqnarray}
for $a_{1}, a_{2}\in A(V),\; f\in \Hom (A(W,z),U),\; w\in A(W,z)$. 

Strengthening Proposition \ref{ptopWU} we have:

\bt{tembeding}
Let $W$ be a weak $V$-module, $U$ a vector space and 
$z$ a nonzero complex number. With the above defined 
$A(V)\otimes A(V)$-module structures,
$\Hom(A(W,z),U)$ and $\Omega({\cal{D}}_{P(z)}(W,U))$ coincide.
\et

\pf Let $f\in \Hom(A(W,z),U)$ and let $v\in V$ be homogeneous.
{}From Proposition \ref{ptopWU}, we have
$$x^{\wt v}Y^{L}(v,x)f,\;\; 
x^{\wt v}Y^{R}(v,x)f\in {\cal{D}}_{P(z)}(W,U)[[x]],$$
Then by expanding $(-z+x)^{\wt v}$ and $(z+x)^{\wt v}$ we get
\begin{eqnarray}
& &\Res_{x}x^{\wt v-1}Y^{R}(v,x)f
=\Res_{x}(-z)^{-\wt v}x^{\wt v-1}(-z+x)^{\wt v}Y^{R}(v,x)f,\\
& &\Res_{x}x^{\wt v-1}Y^{L}(v,x)f
=\Res_{x}z^{-\wt v}x^{\wt v-1}(z+x)^{\wt v}Y^{L}(v,x)f.
\end{eqnarray}
Then for $w\in W$, using (\ref{e3.23}) we have
\begin{eqnarray}
& &\Res_{x}x^{\wt v-1}(Y^{R}(v,x)f)(w)\nonumber\\
&=&\Res_{x}(-z)^{-\wt v}x^{\wt v-1}(-z+x)^{\wt v}
(Y^{R}(v,x)f)(w)\nonumber\\
&=&\Res_{x}(-z)^{-\wt v}x^{\wt v-1}(x-z)^{\wt v}
fY^{o}(v,x)w\nonumber\\
&=&\Res_{x}(-z)^{-\wt v}x^{\wt v-1}(x-z)^{\wt v}
fY(e^{xL(1)}(-x^{-2})^{L(0)}v,x^{-1})w\nonumber\\
&=&\Res_{x}(-1)^{\wt v}(-z)^{-\wt v}x^{-\wt v-1}(x-z)^{\wt v}
fY(e^{xL(1)}v,x^{-1})w\nonumber\\
&=&\Res_{x}(-1)^{\wt v}(-z)^{-\wt v}x^{-1}(1-zx)^{\wt v}
fY(e^{x^{-1}L(1)}v,x)w\nonumber\\
&=&\sum_{i\ge 0}{1\over i!}\Res_{x}(-1)^{\wt v}(-z)^{-\wt v}
x^{-1-i}(1-zx)^{\wt (L(1)^{i}v)+i}
fY(L(1)^{i}v,x)w\nonumber\\
&=&\sum_{i,j\ge 0}{1\over i!}{i\choose j}\Res_{x}(-1)^{\wt v}(-z)^{-\wt v+j}
x^{-1-i+j}(1-zx)^{\wt (L(1)^{i}v)}
fY(L(1)^{i}v,x)w\nonumber\\
&=&\sum_{i\ge 0}{1\over i!}\Res_{x}(-z)^{-\wt v+i}
x^{-1}(1-zx)^{\wt (L(1)^{i}v)}
fY(L(1)^{i}(-1)^{L(0)}v,x)w\nonumber\\
&=&f(\theta(v)*_{P(z)}w).
\end{eqnarray}
Here we are using the fact:
$$\Res_{x}x^{-1-r}(1-zx)^{\wt (L(1)^{i}v)}Y(L(1)^{i}v,x)w\in O(W,z)$$
for $r\ge 1$ (Lemma \ref{lO(W)}).

Similarly, using (\ref{e3.24}) we get
\begin{eqnarray}
& &\Res_{x}x^{\wt v-1}(Y^{L}(v,x)f)(w)\nonumber\\
&=&\Res_{x}z^{-\wt v}x^{\wt v-1}(z+x)^{\wt v}
(Y^{L}(v,x)f)(w)\nonumber\\
&=&\Res_{x}z^{-\wt v}x^{\wt v-1}(x+z)^{\wt v}
fY^{o}(v,x+z)w\nonumber\\
&=&\Res_{x}z^{-\wt v}(x-z)^{\wt v-1}x^{\wt v}
fY^{o}(v,x)w\nonumber\\
&=&\Res_{x}z^{-\wt v}(x-z)^{\wt v-1}x^{\wt v}
fY(e^{xL(1)}(-x^{-2})^{L(0)}v,x^{-1})w\nonumber\\
&=&\Res_{x}(-z)^{-\wt v}(1-zx)^{\wt v-1}x^{-1}
fY(e^{x^{-1}L(1)}v,x)w\nonumber\\
&=&\sum_{i\ge 0}{1\over i!}\Res_{x}(-z)^{-\wt v}x^{-1-i}
(1-zx)^{\wt (L(1)^{i}v)+i-1}
fY(L(1)^{i}v,x)w\nonumber\\
&=&\Res_{x}(-z)^{-\wt v}x^{-1}(1-zx)^{\wt v-1}fY(v,x)w\nonumber\\
& &+\sum_{i\ge 1}{1\over i!}\Res_{x}(-z)^{-\wt v}x^{-2-(i-1)}
(1-zx)^{\wt (L(1)^{i}v)+i-1}
fY(L(1)^{i}v,x)w\nonumber\\
&=&\Res_{x}(-z)^{-\wt v}x^{-1}(1-zx)^{\wt v-1}
fY(v,x)w\nonumber\\
&=&f(w*_{P(z)}v)
\end{eqnarray}
because for $i\ge 1$,
$$\Res_{x}(-1)^{\wt v}x^{-2-(i-1)}(1-zx)^{\wt (L(1)^{i}v)+i-1}
Y(L(1)^{i}v,x)w\in O(W,z).$$
Then it follows immediately from the definitions of the module structures.
$\;\;\;\;\Box$

Recall that $\Hom(A(W,z),U)$ is an $A(V)\otimes A(V)$-module
with the action defined in (\ref{etensormodule}).
Now, {\em let $U$ be a (left) $A(V)$-module instead of just 
a vector space and let $z=-1$}. Then
$\Hom_{A(V)}(A(V),U)$ is a (left) $A(V)$-submodule
of $\Hom (A(V),U)$ equipped with the first action of $A(V)$
(recall (\ref{etensormodule})), i.e.,
\begin{eqnarray}
(af)(b)=f(ba)
\end{eqnarray}
for $a,b\in A(V),\; f\in \Hom(A(V),U)$.
Furthermore, as an $A(V)$-module,
\begin{eqnarray}
U=\Hom_{A(V)}(A(V),U).
\end{eqnarray}

\bd{dinducedmodule}
{\em Let $U$ be an $A(V)$-module.
We define $\Ind_{A(V)}^{V}U$ to be the
$V$-submodule of $({\cal{D}}_{P(-1)}(W,U),Y^{L}_{P(-1)})$, generated by 
$U$ $(=\Hom_{A(V)}(A(V),U))$.}
\ed

We shall briefly use $\Ind\; U$ for $\Ind_{A(V)}^{V}U$
whenever it is clear from the context.

\bl{linducedmodule}
Let $U$ be a (left) $A(V)$-module. Then 
\begin{eqnarray}
U=\Hom_{A(V)}(A(V),U)\subset \Omega(\Ind\; U)\subset \Hom(A(V),U).
\end{eqnarray}
\el

\pf Because  
\begin{eqnarray}
U=\Hom_{A(V)}(A(V),U)\subset \Hom(A(V),U)
=\Omega_{V\otimes V}({\cal{D}}_{P(-1)}(V,U)),
\end{eqnarray}
and the actions $Y^{L}$ and $Y^{R}$ commute, we have
\begin{eqnarray}
Y^{L}(v,x)U\subset (\Omega_{V} ({\cal{D}}_{P(-1)}(V,U),Y^{R}))[[x]]
\end{eqnarray}
for $v\in V$.
Furthermore, because $U$ generates $\Ind\;U$ under the action $Y^{L}$,
we have
\begin{eqnarray}
\Ind \;U  \subset \Omega_{V}({\cal{D}}_{P(-1)}(V,U),Y^{R}).
\end{eqnarray}
Then using Theorem \ref{tembeding}, we get
\begin{eqnarray}
\Omega(\Ind\;U)  \subset 
\Omega_{V}\left(\Omega_{V}({\cal{D}}_{P(-1)}(V,U),Y^{R}),Y^{L}\right)
=\Omega_{V\otimes V}({\cal{D}}_{P(-1)}(V,U))=\Hom(A(V),U).
\end{eqnarray}
This completes the proof.$\;\;\;\;\Box$

Let $U_{1}$ and $U_{2}$ be $A(V)$-modules and left $\psi$
be an $A(V)$-homomorphism from $U_{1}$ to $U_{2}$.
Then $f$ gives rise to a homomorphism $f^{o}$ from $\Hom(V,U_{1})$
to $\Hom(V,U_{2})$ in the obvious way. Furthermore, it is easy 
to see that the restriction of $f^{o}$ is a 
$V\otimes V$-homomorphism from
${\cal{D}}_{P(-1)}(V,U_{1})$ to ${\cal{D}}_{P(-1)}(V,U_{2})$, which
maps $\Hom(A(V),U_{1})$ to $\Hom(A(V),U_{2})$. 
The restriction of $f^{o}$ to $\Ind\; U_{1}$ is a 
$V$-homomorphism from $\Ind\;U_{1}$ to $\Ind\;U_{2}$.
It is routine to check that the map $\Ind: U\mapsto \Ind\; U$ 
gives rise to a functor from the category of 
$A(V)$-modules to the category of weak $V$-modules.
It is also clear that 
\begin{eqnarray}
\Ind (U_{1}\oplus U_{2})=\Ind\; U_{1}\oplus \Ind \; U_{2}.
\end{eqnarray}

Next we study the structure of the induced module $\Ind\; U$.
First, we prove the following result (cf. Lemma \ref{lbasic5}), 
which is a reformulation of a result of [DLM3]:

\bl{lassoc1}
Let $W$ be a weak $V$-module, $w\in W$. Let $u,v\in V$ and
let $k\in {\Z}$ be such that
\begin{eqnarray}
x^{k}Y(u,x)w\in W[[x]],
\end{eqnarray}
or equivalently,
\begin{eqnarray}
u_{k+m}w=0\;\;\mbox{ for }m\ge 0.
\end{eqnarray}
Then for $p,q\in {\Z}$,
\begin{eqnarray}
u_{p}v_{q}w
=\sum_{i=0}^{n}\sum_{j\ge 0}{p-k\choose i}{k\choose j}
(u_{p-k-i+j}v)_{q+k+i-j}w.
\end{eqnarray}
where $n$ is any nonnegative integer such that $x^{n+1+q}Y(v,x)w\in W[[x]]$. 
\el

\pf Since $x^{k}Y(u,x)w\in W[[x]]$, by applying 
$\Res_{x_{1}}x_{1}^{k}$ to the Jacobi identity for the triple
$(u,v,w)$ we get [DL]
\begin{eqnarray}\label{eweakassoc10}
(x_{0}+x_{2})^{k}Y(u,x_{0}+x_{2})Y(v,x_{2})w=
(x_{2}+x_{0})^{k}Y(Y(u,x_{0})v,x_{2})w.
\end{eqnarray}
Notice that
\begin{eqnarray}
u_{p}v_{q}w
=\Res_{x_{0}}\Res_{x_{2}}x_{2}^{q}(x_{0}+x_{2})^{p}
Y(u,x_{0}+x_{2})Y(v,x_{2})w.
\end{eqnarray}
We can multiply the left hand side of (\ref{eweakassoc10}) by
$x_{2}^{q}(x_{0}+x_{2})^{p-k}$,
but we are not allowed to multiply the right hand side of 
(\ref{eweakassoc10}) by
$x_{2}^{q}(x_{0}+x_{2})^{p-k}$. Notice that
\begin{eqnarray}
(x_{0}+x_{2})^{p-k}
=\sum_{i=0}^{n}{p-k\choose i}x_{0}^{p-k-i}x_{2}^{i}
+\sum_{i\ge n+1}{p-k\choose i}x_{0}^{p-k-i}x_{2}^{i}
\end{eqnarray}
and
\begin{eqnarray}\label{e713t}
\Res_{x_{2}}\sum_{i\ge n+1}{p-k\choose i}x_{0}^{p-k-i}x_{2}^{i+q}
(x_{0}+x_{2})^{k}Y(u,x_{0}+x_{2})Y(v,x_{2})w=0
\end{eqnarray}
because $x^{n+1+q}Y(v,x)w\in W[[x]]$.
Using (\ref{eweakassoc10})-(\ref{e713t}) we get
\begin{eqnarray}
& & u_{p}v_{q}w\nonumber\\
&=&\Res_{x_{0}}\Res_{x_{2}}(x_{0}+x_{2})^{p}x_{2}^{q}
Y(u,x_{0}+x_{2})Y(v,x_{2})w\nonumber\\
&=&\Res_{x_{0}}\Res_{x_{2}}(x_{0}+x_{2})^{p-k}x_{2}^{q}
\left((x_{0}+x_{2})^{k}Y(u,x_{0}+x_{2})Y(v,x_{2})w\right)\nonumber\\
&=&\Res_{x_{0}}\Res_{x_{2}}\sum_{i=0}^{n}{p-k\choose i}
x_{0}^{p-k-i}x_{2}^{i+q}
\left((x_{0}+x_{2})^{k}Y(u,x_{0}+x_{2})Y(v,x_{2})w\right)\nonumber\\
&=&\Res_{x_{0}}\Res_{x_{2}}\sum_{i=0}^{n}{p-k\choose i}
x_{0}^{p-k-i}x_{2}^{i+q}
\left((x_{2}+x_{0})^{k}Y(Y(u,x_{0})v,x_{2})w\right)\nonumber\\
&=&\sum_{i=0}^{n}\sum_{j\ge 0}{p-k\choose i}{k\choose j}
(u_{p-k-i+j}v)_{q+k+i-j}w.
\end{eqnarray}
This concludes the proof. $\;\;\;\;\Box$

As an immediate consequence of Lemma \ref{lassoc1}, 
we have the following result,
which was proved in [DM] and [Li1].

\bc{cdmlif}
Let $W$ be a weak $V$-module and let $w\in W$. Set
\begin{eqnarray}
\<w\>=\mbox{{\rm linear span }}\{ v_{m}w\;|\; v\in V,\; m\in {\Z}\}.
\end{eqnarray}
Then $\<w\>$ is the sub-weak-module of $W$, generated by $w$.
\ec

Furthermore, we have:

\bl{lgenerating}
Let $W$ be a weak $V$-module and let $U$ be an irreducible 
$A(V)$-submodule of $\Omega (W)$.  
Then the weak submodule $\<U\>$ of $W$, generated by $U$,
is a lowest weight generalized $V$-module with $U$
being the lowest weight subspace.
\el

\pf From Corollary \ref{cdmlif} we have
\begin{eqnarray}\label{espanf}
\<U\>=\mbox{linear span} \{ v_{m}U\;|\; v\in V,\; m\in {\Z}\}.
\end{eqnarray}
Since $V$ has a countable basis, $A(V)$, being a quotient of $V$,
also has a countable basis. Then 
the central element $\omega+O(V)$
of $A(V)$ acts on any irreducible $A(V)$-module as a scalar.
(See the proof of Lemma 1.2.1, [Z].) That is, $L(0)$ 
acts as a scalar $h$ on $U$ (being a subspace of $W$). Then 
it immediately follows from (\ref{espanf}) and the following facts:
\begin{eqnarray}
& &\wt v_{m}=\wt v-m-1,\\
& &v_{n}U=0
\end{eqnarray}
for homogeneous $v\in V$ and for $m\in {\Z},\;n\ge \wt v$.$\;\;\;\;\Box$

We immediately have:

\bc{cL(U)10}
Let $U$ be an irreducible (left) $A(V)$-module. Then
$\Ind\; U$ is a lowest weight generalized $V$-module
with $U$ as the lowest weight subspace. $\;\;\;\;\Box$
\ec

\br{rL(U)2}
{\em There are two questions regarding $\Ind\; U$:
(1) Is $\Ind_{A(V)}^{V}U$ already an irreducible generalized $V$-module?
(2) Is there a canonical characterization of $\Ind\;U$?}
\er

\section{Functor $F$ and Frenkel-Zhu's fusion rule theorem}
The main goal of this section is to give an alternate proof
of Frenkel and Zhu's fusion rule theorem.

Recall from [B] (cf. [FFR], [Li1]) the Lie algebra $g(V)$ associated to
the vertex operator algebra $V$. As a vector space,
$$g(V)=\hat{V}/D\hat{V},$$
where 
$$\hat{V}=V\otimes {\C}[t,t^{-1}],\;\;\;
D=L(-1)\otimes 1+1\otimes {d\over dt}.$$
The Lie bracket is given by
$$[u(m),v(n)]=\sum_{i\ge 0}{m\choose i}(u_{i}v)(m+n-i)$$
for $u,v\in V,\; m, n\in {\Z}$, where $u(m)=u\otimes t^{m}$. 
Furthermore, $g(V)$ is naturally a ${\Z}$-graded Lie algebra with
\begin{eqnarray}
\deg v(m)=\wt v-m-1
\end{eqnarray}
for homogeneous $v\in V$ and for $m\in {\Z}$.
It is clear that any weak $V$-module
is a natural $g(V)$-module and that any generalized $V$-module is a 
${\C}$-graded $g(V)$-module. 
It was known (cf. [Li1]) that $A(V)_{Lie}$ is a natural quotient 
Lie algebra of $g(V)_{0}$, where $g(V)_{0}$ is the degree-zero 
Lie subalgebra of $g(V)$. Then any $A(V)$-module is 
a natural $g(V)_{0}$-module.

Recall a notion from [DLM2].
(Here we use a different symbol for the universal object.)

\bd{duniversalmodule}
{\em Let $U$ be an $A(V)$-module. Then $U$ is a $g(V)_{0}$-module.
View $U$ as a $(g(V)_{0}+g(V)_{-})$-module with $g(V)_{-}U=0$, where
$g(V)_{-}=\oplus_{n>0}g(V)_{-n}$.
Define the standard induced $g(V)$-module
\begin{eqnarray}
\tilde{F}(U)=U(g(V))\otimes_{U(g(V)_{-}+g(V)_{0})}U,
\end{eqnarray}
which is an ${\N}$-graded $g(V)$-module with
\begin{eqnarray}
\deg U=0.
\end{eqnarray}
Then we define $F(U)$ to be the quotient $g(V)$-module of $\tilde{F}(U)$
modulo the following Jacobi identity relation:
\begin{eqnarray}\label{ejacobirelation}
& &x_{0}^{-1}\delta\left(\frac{x_{1}-x_{2}}{x_{0}}\right)
Y(u,x_{1})Y(v,x_{2})w
-x_{0}^{-1}\delta\left(\frac{x_{2}-x_{1}}{-x_{0}}\right)
Y(v,x_{2})Y(u,x_{1})w\nonumber\\
&=&x_{2}^{-1}\delta\left(\frac{x_{1}-x_{0}}{x_{2}}\right)
Y(Y(u,x_{0})v,x_{2})w
\end{eqnarray}
for $u,v\in V,\;w\in \tilde{F}(U)$.}
\ed

{}From definition, $F(U)$ is an ${\N}$-graded $g(V)$-module.
Because of (\ref{ejacobirelation}), $F(U)$ clearly
is an ${\N}$-graded weak $V$-module.
Let $e_{U}$ be the natural map from $U$ to $F(U)$. Then
we have the following obvious universal property: 

\bp{puniversal}
Let $W$ be any weak $V$-module and let $\psi$ be 
any $A(V)$-homomorphism from $U$ to $\Omega(W)$. Then
there exists a unique $V$-homomorphism $\tilde{\psi}$ from 
$F(U)$ to $W$ such that $\tilde{\psi}e_{U}=\psi$. $\;\;\;\;\Box$
\ep

Note that we have not excluded the possibility that $F(U)=0$ 
even if $U\ne 0$. With the weak $V$-module $\Ind \;U$ we have
the following result:

\bl{ltograded}
Let $U$ be an $A(V)$-module. Then the natural linear map 
$e_{U}$ from $U$ to $F(U)$ is injective and $e_{U}(U)=F(U)(0)$.
\el

\pf It is clear that $e_{U}(U)=F(U)(0)$. 
Since $U$ is an $A(V)$-submodule of $\Omega(\Ind\;U)$,
using the universal property of $F(U)$ (Proposition \ref{puniversal}),
we obtain a 
$V$-homomorphism $\phi$ from $F(U)$ to $\Ind\;U$ such that
$\phi e_{U}$ is the embedding of $U$ into $\Omega(\Ind\;U)$.
In particular, $\phi e_{U}$ is injective. Consequently,
$e_{U}$ is injective. $\;\;\;\;\Box$

In view of Lemma \ref{ltograded}, we consider $U$ as a 
canonical subspace of $F(U)$.
Combining Lemma \ref{ltograded} with Lemma \ref{lgenerating} 
we immediately have:

\bl{lfree2}
Let $U$ be an irreducible $A(V)$-module. Then $F(U)$ is a lowest weight
generalized $V$-module with $U$ as the lowest weight subspace. $\;\;\;\;\Box$
\el

Consider all graded submodules $W$ of $F(U)$ such that $W\cap U=0$.
Then the sum of all such graded submodules is still a such graded submodule,
so that it is the unique maximal graded submodule with this property.
Define $L(U)$ to be the quotient module of $F(U)$ modulo 
the maximal submodule. We have ([Z], Theorem 2.2.1):

\bl{lL(U)}
Let $U$ be an $A(V)$-module. Then $L(U)$ is an ${\N}$-graded 
weak $V$-module such that for any
nonzero graded submodule $W$ of $L(U)$, $U\cap W\ne 0$.
Furthermore, if $U$ is irreducible, $L(U)$ is an
irreducible generalized $V$-module.
\el

\pf The first assertion directly follows from the definition of $L(U)$.
Since $U$ is irreducible, by Lemma \ref{lgenerating}, 
$L(U)$ is a lowest weight 
generalized $V$-module with $U$ as the lowest weight subspace.
Then the ${\N}$-grading on $L(U)$ is a shift of the $L(0)$-grading
on $W$. Consequently,
any submodule of $L(U)$ is automatically graded. It follows immediately that
$L(U)$ is irreducible. $\;\;\;\;\Box$

It is routine to check that the map $F: U\mapsto F(U)$ 
gives rise to a functor $F$ from the category of $A(V)$-modules 
to the category of ${\N}$-graded weak $V$-modules.
Furthermore, given a family of $A(V)$-modules $U_{i}$ for $i\in S$, 
we have
\begin{eqnarray}
F(\oplus_{i\in S}U_{i})=\oplus_{i\in S}F(U_{i}),
\end{eqnarray}
or equivalently, if $U$ is an $A(V)$-module such that
$U=E\otimes U_{1}$ where $E$ is a vector space and $U_{1}$
is an $A(V)$-module, then
$F(U)=E\otimes F(U_{1})$. We also have the following analogue of the
Frobenius reciprocity theorem (cf. [Ki]):

\bl{lFrobeniusF(U)}
Let $W$ be a weak $V$-module and let $U$ be an $A(V)$-module. Then
the map 
\begin{eqnarray}
\Omega':& & \Hom_{V}(F(U),W)\rightarrow \Hom_{A(V)}(U,\Omega(W))\nonumber\\
& &\psi\mapsto \Omega(\psi)
\end{eqnarray}
is a linear isomorphism.
\el

\pf Because $U$ generates $F(U)$ as a weak $V$-module, it is clear that
$\Omega'$ is injective. It follows from the universal property of
$F(U)$ (Proposition \ref{puniversal}) that $\Omega'$ is also surjective.
$\;\;\;\;\Box$

\br{rFrobeniusIndU}
{\em Let $W$ and $U$ be given as in Lemma \ref{lFrobeniusF(U)}.
Similarly, we define a linear map $\Omega''$ from $\Hom_{V}(\Ind\;U,W)$
to $\Hom_{A(V)}(U,\Omega(W))$. Then $\Omega''$ is injective.
It is easy to see that $\Omega''$ is surjective if and only if
the $V$-homomorphism from $F(U)$ to $\Ind\;U$, extending the identity
map of $U$, is an isomorphism.}
\er

We shall need the following fact:

\bl{ltensoruniversal}
Let $V_{1}$ and $V_{2}$ be vertex operator algebras and let $U_{1}$
and $U_{2}$ be $A(V_{1})$ and $A(V_{2})$-modules, respectively.
Let $W$ be a weak $V_{1}\otimes V_{2}$-module and let $\psi$ be an
$A(V_{1})\otimes A(V_{2})$-homomorphism from $U_{1}\otimes U_{2}$
to $\Omega(W)$. Then there exists a unique 
$V_{1}\otimes V_{2}$-homomorphism $\bar{\psi}$ from $F(U_{1})\otimes F(U_{2})$ 
to $W$, extending $\psi$.
\el

\pf The uniqueness is clear because $U_{1}\otimes U_{2}$ generates
$F(U_{1})\otimes F(U_{2})$ as a weak $V_{1}\otimes V_{2}$-module.
By Proposition \ref{puniversal}, there exists a $V_{1}$-homomorphism 
$\psi_{1}$ from $F_{V_{1}}(U_{1}\otimes U_{2})$ to $W$, extending $\psi$. 
Note that
$$F_{V_{1}}(U_{1})\otimes U_{2}=F_{V_{1}}(U_{1}\otimes U_{2}).$$ 
It is clear that $\psi_{1}$ is an $A(V_{2})$-homomorphism.
Then by Proposition \ref{puniversal} again, there exists a 
$V_{2}$-homomorphism $\psi_{2}$ from 
$$F_{V_{2}}(F(U_{1})\otimes U_{2})
\;\left(=F_{V_{1}}(U_{1})\otimes F_{V_{2}}(U_{2})\right)$$ to $W$,
extending $\psi_{1}$. Consequently, $\psi_{2}$ is  a 
$V_{1}\otimes V_{2}$-homomorphism $\bar{\psi}$ from $F(U_{1})\otimes F(U_{2})$ 
to $W$, extending $\psi$. $\;\;\;\;\Box$

\br{rtensoruniversal} {\em 
Let $V_{1}, V_{2}, U_{1}$ and $U_{2}$ be given as in Lemma \ref{ltensoruniversal}.
It was proved in [DMZ] that $A(V_{1}\otimes V_{2})$ is naturally isomorphic
to $A(V_{1})\otimes A(V_{2})$. Then
$U_{1}\otimes U_{2}$ is a natural $A(V_{1}\otimes V_{2})$-module. 
By Lemma \ref{ltensoruniversal}, there exists a 
$V_{1}\otimes V_{2}$-homomorphism $\psi$ from $F(U_{1})\otimes F(U_{2})$
to $F_{V_{1}\otimes V_{2}}(U_{1}\otimes U_{2})$, extending
the identity map of $U_{1}\otimes U_{2}$. It follows from
the universal property of $F_{V_{1}\otimes V_{2}}(U_{1}\otimes U_{2})$
that $\psi$ is an isomorphism.}
\er

Recall the involution (anti-automorphism)  $\theta$ of $A(V)$.
Let $U$ be a (left) $A(V)$-module. 
Then from the classical fact $U^{*}$ is a left $A(V)$-module 
with the action defined by
\begin{eqnarray}
(af)(u)=f(\theta(a)u)\;\;\;\mbox{ for }a\in A(V),\; u\in U.
\end{eqnarray}
The following result is classical in nature.

\bl{lUUembedding}
Let $U_{1}$, $U_{2}$ be (left) $A(V)$-modules and let $B$ be an 
$A(V)$-bimodule. Define
\begin{eqnarray}
d: & &\Hom (B\otimes _{A(V)}U_{1}, U_{2})\rightarrow 
\Hom (U_{1}\otimes U_{2}^{*}, B^{*})\nonumber\\
& &\psi\mapsto d_{\psi},
\end{eqnarray}
where for $u_{1}\in U_{1},\; u_{2}^{*}\in U_{2}^{*}$,
\begin{eqnarray}
\<d_{\psi}(u_{1}\otimes u_{2}^{*}),b\>
=\<u_{2}^{*},\psi(b\otimes u_{1})\>.
\end{eqnarray}
Then 
\begin{eqnarray}
d\left(\Hom_{A(V)} (B\otimes _{A(V)}U_{1}, U_{2})\right)\subset 
\Hom_{A(V)\otimes A(V)} (U_{1}\otimes U_{2}^{*}, B^{*}),
\end{eqnarray}
where $B^{*}$ is considered as an $A(V)\otimes A(V)$-module
with the action defined by (\ref{etensormodule}).
If we in addition assume that $U_{2}$ is finite-dimensional,
then the restriction of $d$ gives rise to a linear isomorphism from
$\Hom_{A(V)}(B\otimes _{A(V)}U_{1}, U_{2})$ onto 
$\Hom_{A(V)\otimes A(V)}(U_{1}\otimes U_{2}^{*}, B^{*})$.

In particular, let $U$ be a (left) $A(V)$-module. Then the map
\begin{eqnarray}
d_{U}: & &U\otimes U^{*}\rightarrow A(V)^{*}\nonumber\\
& &u\otimes u^{*}\mapsto d_{U}(u\otimes u^{*}),
\end{eqnarray}
where for $a\in A(V)$,
\begin{eqnarray}
\<d_{U}(u\otimes u^{*}),a\>=\<u^{*},au\>,
\end{eqnarray}
is an $A(V)\otimes A(V)$-homomorphism.
\el

\pf Let $\eta'$ be the natural embedding of $\Hom (B\otimes U_{1}, U_{2})$
into $\Hom (U_{1}\otimes U_{2}^{*}, B^{*})$. It is a classical fact that
$\eta'$ is a linear isomorphism if $U_{2}$ is finite-dimensional.
With $B\otimes_{A(V)} U_{1}$
being a quotient space of $B\otimes U_{1}$, we naturally consider
$\Hom (B\otimes_{A(V)} U_{1}, U_{2})$ as a subspace of 
$\Hom (B\otimes U_{1}, U_{2})$. Then $d$ is the restriction of $\eta'$.
Thus $d$ is injective. 
Let $\psi\in \Hom(B\otimes U_{1}, U_{2}),\;
a_{1}, a_{2}\in A(V)$, $u_{1}\in U_{1}, 
\;u_{2}^{*}\in U_{2}^{*}$ and $b\in B$. Then
\begin{eqnarray}
\<d_{\psi}((a_{1},a_{2})(u_{1}\otimes u_{2}^{*})),b\>
&=&\<d_{\psi}(a_{1}u_{1}\otimes a_{2}u^{*}_{2}),b\>\nonumber\\
&=&\<a_{2}u^{*}_{2},\psi(b\otimes a_{1}u_{1})\>\nonumber\\
&=&\<u^{*}_{2},\theta(a_{2})\psi(b\otimes a_{1}u_{1})\>
\end{eqnarray}
and
\begin{eqnarray}
& &\<(a_{1},a_{2})d_{\psi}(u_{1}\otimes u^{*}_{2}),b\>\>\nonumber\\
&=&\<d_{\psi}(u_{1}\otimes u^{*}_{2}),\theta(a_{2})ba_{1}\>\nonumber\\
&=&\<u^{*}_{2},\psi(\theta(a_{2})ba_{1}\otimes u_{1})\>.
\end{eqnarray}
It follows immediately that
$\psi\in \Hom_{A(V)} (B\otimes_{A(V)} U_{1}, U_{2})$ if and only if
$$d_{\psi}\in \Hom_{A(V)\otimes A(V)}(U_{1}\otimes U_{2}^{*}, B^{*}).$$
This completes the proof.$\;\;\;\;\Box$

\br{rinducedmodule}
{\em Let $U$ be an irreducible (left) $A(V)$-module. 
Let $u^{*}$ be a nonzero element of $U^{*}$.
By Lemma \ref{lUUembedding},
$d_{U}(\cdot \otimes u^{*})$ gives an $A(V)$-homomorphism from $U$ to 
$A(V)^{*}$ equipped with the first action. It follows from the 
irreducibility of $U$
that $d_{U}(\cdot \otimes u^{*})$ is injective.
The it follows from Lemma \ref{lgenerating} that there is a canonical lowest 
weight generalized $V$-module
inside the regular representation on ${\cal{D}}_{P(-1)}(V)$
with $U$ as the lowest weight subspace.}
\er

\bl{lintertwiningoperator}
Let $W$ be a weak $V$-module and let $W_{1}$ and $W_{2}$
be lowest weight generalized $V$-modules with lowest weight subspaces 
$W_{1}(0)$ and $W_{2}(0)$, respectively.
Define the restriction map
\begin{eqnarray}
\Omega_{T}: & &\Hom_{V\otimes V}(W_{1}\otimes W_{2}, {\cal{D}}_{P(-1)}(W))
\rightarrow 
\Hom_{A(V)\otimes A(V)}(W_{1}(0)\otimes W_{2}(0),
 \Omega({\cal{D}}_{P(-1)}(W)))\nonumber\\
& & \psi \mapsto \Omega(\psi)|_{W_{1}(0)\otimes W_{2}(0)}.
\end{eqnarray}
Then $\Omega_{T}$ is injective.
\el

\pf It immediately follows from the fact that
the $V\otimes V$-module $W_{1}\otimes W_{2}$ is generated by
$W_{1}(0)\otimes W_{2}(0)$. $\;\;\;\;\Box$

Recall from Theorem \ref{trecall4} that for generalized $V$-modules 
$W,W_{1}$ and $W_{2}$, $F_{p}[P(z)]_{W_{1}W_{2}}^{W'}$ is 
a linear isomorphism from ${\cal{V}}_{W_{1}W_{2}}^{W'}$ onto 
$\Hom_{V\otimes V}(W_{1}\otimes W_{2}, {\cal{D}}_{P(z)}(W))$.
For the rest of this section, we use $F_{W_{1}W_{2}}^{W'}$ for
$F_{0}[P(-1)]_{W_{1}W_{2}}^{W'}$.

Combining Theorem \ref{trecall4} and Lemma \ref{lintertwiningoperator} 
with Theorem \ref{tembeding} we immediately have:

\bp{pinequality}
Let $W, W_{1}$ and $W_{2}$ be lowest weight generalized $V$-modules.
Then 
$\Omega_{T} F_{W_{1}W_{2}}^{W'}$
is an injective linear map from ${\cal{V}}^{W'}_{W_{1}W_{2}}$
to $\Hom_{A(V)\otimes A(V)}(W_{1}(0)\otimes W_{2}(0), A(W)^{*})$.
In particular,
\begin{eqnarray}
\dim {\cal{V}}^{W'}_{W_{1}W_{2}}\le 
\dim \Hom_{A(V)\otimes A(V)}(W_{1}(0)\otimes W_{2}(0), A(W)^{*}).
\end{eqnarray}
\ep

Next, we shall show that $\Omega_{T} F_{W_{1}W_{2}}^{W'}$ 
is a linear isomorphism in a certain situation.

\bt{tinducedintertwiningoperator}
Let $W$ be a lowest weight generalized $V$-module 
and let $U_{1},\; U_{2}$ be 
finite-dimensional irreducible $A(V)$-modules. 
Then the linear map 
$$\Omega_{T} F_{F(U_{1})F(U_{2})}^{W'}:\;\;
 {\cal{V}}^{W'}_{F(U_{1}) F(U_{2})}\rightarrow 
\Hom_{A(V)\otimes A(V)}(U_{1}\otimes U_{2},A(W)^{*})$$
is a linear isomorphism.
\et

\pf We only need to prove that $\Omega_{T} F_{F(U_{1})F(U_{2})}^{W'}$ 
is onto. For simplicity, in this proof we use ${\cal{F}}$ for
$F_{F(U_{1})F(U_{2})}^{W'}$. Let
$$\psi\in \Hom_{A(V)\otimes A(V)}(U_{1}\otimes U_{2},A(W)^{*}).$$
Then $\psi(U_{1}\otimes U_{2})$ is an 
$A(V)\otimes A(V)$-submodule of $A(W)^{*}$, which is
$\Omega({\cal{D}}_{P(-1)}(W))$ by Proposition \ref{ptopWU} 
with $U={\C}$ and $z=-1$.
By Lemma \ref{lfree2} we may assume that 
$\omega+O(V)$ acts as scalars $h_{1}$ and $h_{2}$ on
$U_{1}$ and $U_{2}$, respectively. Then $L(0)$ acts as scalar
$h_{1}+h_{2}$ on $\psi(U_{1}\otimes U_{2})$.
Let $E$ be the $V\otimes V$-submodule of ${\cal{D}}_{P(-1)}(W)$, generated by
$\psi(U_{1}\otimes U_{2})$. 
By Lemma \ref{ltensoruniversal},
$\psi$ extends to a $V\otimes V$-homomorphism $\bar{\psi}$ from
$F(U_{1})\otimes F(U_{2})$ to $E$. By Theorem  \ref{trecall4}, 
we get an intertwining
operator ${\cal{F}}^{-1}(\bar{\psi})$ of type 
${W'\choose F(U_{1}) F(U_{2})}$ such that
$$\Omega_{T} {\cal{F}}({\cal{F}}^{-1}(\bar{\psi})=\Omega(\bar{\psi})=\psi.$$
Thus $\Omega_{T} {\cal{F}}$ is onto. This completes the proof.$\;\;\;\;\Box$

\br{rT} 
{\em In the proof, {\em if} $E$ is an irreducible generalized
$V\otimes V$-module, 
then from [FHL], $E=L(U_{1})\otimes L(U_{2})$, and then
$\Omega {\cal{F}}$ will be a linear isomorphism from 
${\cal{V}}^{W'}_{L(U_{1}) L(U_{2})}$ to 
$\Hom_{A(V)\otimes A(V)}(U_{1}\otimes U_{2},A(W)^{*})$.
But, $E$ in general is not irreducible.}
\er

Combining Theorem \ref{tinducedintertwiningoperator}
with  Lemma \ref{lUUembedding} we immediately have:

\bc{cinducedintertwiningoperator}
Let $W, U_{1}$ and $U_{2}$ be as in 
Theorem \ref{tinducedintertwiningoperator}.
Then $d^{-1}\Omega_{T} F_{F(U_{1})F(U_{2}^{*})}^{W'}$ 
is a linear isomorphism from 
${\cal{V}}^{W'}_{F(U_{1}) F(U_{2}^{*})}$
to $\Hom_{A(V)\otimes A(V)}(A(W)\otimes_{A(V)}U_{1}, U_{2})$.
$\;\;\;\;\Box$
\ec

Now we immediately have the following modified 
Frenkel-Zhu's fusion rule theorem (cf. [FZ], [Li1] and [Li2]):

\bc{cintertwining}
Let $W$ be a $V$-module and let $W_{1}$ and $W_{2}$ be irreducible
$V$-modules such that $W_{1}=F(W_{1}(0))$ and
$W_{2}'=F((W_{2}(0))^{*})$, or what is equivalent,
$F(W_{1}(0))$ and $F((W_{2}(0))^{*})$
 are irreducible. Then
\begin{eqnarray}
\dim \Hom_{A(V)}(A(W)\otimes_{A(V)}W_{1}(0),W_{2}(0))
=\dim {\cal{V}}^{W_{2}}_{WW_{1}}.
\end{eqnarray}
In particular, this is true if $V$ satisfies the condition that 
every lowest weight generalized $V$-module is completely reducible.
\ec

\pf It was proved in [HL2] (cf. [FHL]) that 
\begin{eqnarray}
\dim {\cal{V}}^{W_{2}}_{WW_{1}}=\dim {\cal{V}}^{W_{2}}_{W_{1}W}
=\dim {\cal{V}}^{W'}_{W_{1}W_{2}'}.
\end{eqnarray}
Then it follows immediately from 
Corollary \ref{cinducedintertwiningoperator}.
$\;\;\;\;\Box$

\end{document}